\documentclass[paper,twocolumn,twoside]{geophysics}

\usepackage{latexsym,amssymb,amsmath,amsthm, amsfonts, float}
\usepackage{algorithmicx}
\usepackage{algpseudocode}

\usepackage[draft]{hyperref}

\usepackage{enumerate}

\usepackage{setspace}

\usepackage{color}

\def\cal{\mathcal}

\def\cH{{\cal H}}

\def\cO{{\cal O}}

\newtheorem{algorithm}{Algorithm}

\newcommand{\y}{\mathbf{y}}
\renewcommand{\v}{\mathbf{v}}

\renewcommand{\u}{\mathbf{u}}
\newcommand{\f}{\mathbf{f}}

\newcommand{\x}{\mathbf{x}}

\newcommand{\beq}{\begin{equation}}
\newcommand{\eeq}{\end{equation}}
\renewcommand{\tilde}{\widetilde}
\renewcommand{\hat}{\widehat}
\newcommand{\bit}{\begin{itemize}}
\newcommand{\eit}{\end{itemize}}
\newcommand{\ben}{\begin{enumerate}}
\newcommand{\een}{\end{enumerate}}

\title{The method of polarized traces for the 3D Helmholtz equation}

\address{
\footnotemark[1] Massachusetts Institute of Technology, \\
Department of Mathematics and Earth Resources Laboratory, \\
77 Massachusetts Ave, \\
Cambridge, MA 02139. \\
\footnotemark[2] Universit\'e Catholique de Louvain, \\
1, Place de l'Universit\'e, \\
B-1348 Louvain-la-Neuve, Belgium\\
\footnotemark[3] Total E. \& P. Research \& Technology USA\\
1201 Louisiana St.,\\
Houston, TX 77002.\\
\footnotemark[4] Computational Research Division, \\
Lawrence Berkeley National Laboratory,\\
1 Cyclotron road, \\
Berkeley, CA 94720. \\
\footnotemark[5] University of California Irvine, \\
Department of Mathematics,\\
540 Rowland Hall, \\
Irvine, CA 92963.}
\footer{Example}
\lefthead{Zepeda-N\'u\~nez et al.}
\righthead{Polarized Traces for 3D Helmholtz}

\author{ Leonardo Zepeda-N\'u\~nez\footnotemark[1]\footnotemark[4]\footnotemark[5], Adrien Scheuer\footnotemark[1]\footnotemark[2], Russell J. Hewett\footnotemark[3], and Laurent Demanet\footnotemark[1]}

\date{December 2017}

\begin{document}

\maketitle

\begin{abstract}
We present a fast solver for the 3D high-frequency Helmholtz equation in heterogeneous, constant density, acoustic media. The solver is based on the method of polarized traces, coupled with distributed linear algebra libraries and pipelining to obtain an empirical online runtime $ \cO(\max(1,R/n) N \log N)$ where $N = n^3$ is the total number of degrees of freedom and $R$ is the number of right-hand sides. Such a favorable scaling is a prerequisite for large-scale implementations of full waveform inversion (FWI) in frequency domain.

\end{abstract}

\section{Introduction} \label{section:introduction}

Efficient modeling of time-harmonic wave scattering in heterogeneous acoustic or elastic media remains a difficult problem in numerical analysis, yet it has broad application in seismic inversion techniques, as shown by \cite{Chen:Inverse_scattering_via_Heisenberg's_uncertainty_principle,Pratt:Seismic_waveform_inversion_in_the_frequency_domain;_Part_1_Theory_and_verification_in_a_physical_scale_model,Virieux_Operto:An_overview_of_full-waveform_inversion_in_exploration_geophysics}.
In the constant density acoustic approximation, time-harmonic wave propagation is modeled by the Helmholtz equation,
\begin{equation}
    \triangle u(\x) + \omega^2 m (\x) u(\x) = f_s(\x), \qquad \mbox{in } \Omega \label{eq:Helmholtz}
\end{equation}
with absorbing boundary conditions,  and where  $\Omega $ is a $3$D rectangular domain, $\triangle$ is the 3D Laplacian, $\x = (x,y,z)$, $m = 1/c^2(\x)$ is the squared slowness for velocity $c(\x)$, $u$ is the wavefield, and $f_s$ are the sources, indexed by $s = 1,..., R$. There is no essential obstruction to extending the techniques presented in this paper to the case of heterogeneous density, or to the elastic, viscoelastic or poroelastic time-harmonic equations.

Throughout this paper we assume that Eq.~\ref{eq:Helmholtz} is in the high-frequency regime, i.e., when $\omega \sim n$, where $n$ is the number of unknowns in each dimension.  Alternatively, in the situation where the domain $\Omega$ grows and the frequency is fixed, the problem can be rescaled in such a way that $\omega$ grows in proportion to the domain size, which can also considered to be a ``high-frequency'' regime.  For this paper, we focus on the former scenario.

Given the importance of solving Eq.~\ref{eq:Helmholtz} in geophysical contexts, there has been a renewed interest in developing efficient algorithms to solve the ill-conditioned linear system resulting from its discretization.   Recent progress toward an efficient solver, i.e., a solver with linear complexity, has generally focused on three strategies:
\begin{itemize}
\item \textit{Fast direct solvers}, such as the ones introduced by \cite{Xia:multifrontal,Wang:H_multifrontal,Gillman_Barnett_Martinsson:A_spectrally_accurate_solution_technique_for_frequency_domain_scattering_problems_with_variable_media,Amestoy_Weisbecker:compressed_MUMPS}, which couple multifrontal techniques (e.g., \cite{GeorgeNested_dissection,Duff_Reid:The_Multifrontal_Solution_of_Indefinite_Sparse_Symmetric_Linear}) with compressed linear algebra (e.g., \cite{Bebendorf:2008}) to obtain efficient direct solvers with small memory footprint.  However, they suffer the same sub-optimal asymptotic complexity as standard multifrontal methods (e.g., \cite{Demmel_Li:superlu,Amestoy_Duff:MUMPS,Davis:UMFPACK}) in the high-frequency regime.

\item \textit{Classical preconditioners}, such as {incomplete factorization preconditioners} (e.g., \cite{Grote_Schenk:algebraic_multilever_preconditioner_Helmholtz_equation}) and {multigrid-based preconditioners} (e.g., \cite{Brandt_Livshits:multi_ray_multigrid_standing_wave_equations,Erlangga:shifted_laplacian,Sheikh_Lahaye_Vuik:On_the_convergence_of_shifted_Laplace_preconditioner_combined_with_multilevel_deflation,Calandra_Grattonn:an_improved_two_grid_preconditioner_for_the_solution_of_3d_Helmholtz}), which are relatively simple to implement but suffer from super-linear asymptotic complexity and may need significant tuning to achieve effective run-times.

\item \textit{Sweeping-like preconditioners} (e.g., \cite{GanderNataf:LU_incomplete,EngquistYing:Sweeping_H,EngquistYing:Sweeping_PML,Chen_Xiang:a_source_transfer_ddm_for_helmholtz_equations_in_unbounded_domain,Cheng_Xiang:A_Source_Transfer_Domain_Decomposition_Method_For_Helmholtz_Equations_in_Unbounded_Domain_Part_II_Extensions,CStolk_rapidily_converging_domain_decomposition,GeuzaineVion:double_sweep,Liu_Ying:Recursive_sweeping_preconditioner_for_the_3d_helmholtz_equation,ZepedaDemanet:the_method_of_polarized_traces}), which are a relatively recent domain decomposition based approach that has been shown to achieve linear or nearly-linear asymptotic complexity.
\end{itemize}

The method in this paper belongs to the third category. Sweeping preconditioners and their generalizations, i.e., domain decomposition techniques coupled with high-quality transmission/absorption conditions, offer the right mix of ideas to attain linear or near-linear complexity in 2D and 3D, provided that the medium does not have large resonant cavities~\citep{ZepedaDemanet:the_method_of_polarized_traces}. These methods rely on the sparsity of the linear system to decompose the domain in layers, in which classical sparse direct methods are used to compute the interactions within the layer. Interactions across layers are computed by sequentially sweeping through the sub-domains in an iterative fashion.

For current applications, empirical runtimes are a more practical measure of an algorithm's performance than asymptotic complexity. This requirement has led to a recent effort to reduce the runtimes of preconditioners with optimal asymptotic complexity by leveraging parallelism. For example, \citet{Poulson_Engquist:a_parallel_sweeping_preconditioner_for_heteregeneous_3d_helmholtz} introduce a new local solver, i.e., a solver for the subproblem defined on each layer, which carefully handles communication patterns between layers to obtain impressive timings.  While most sweeping algorithms require visiting each subdomain in sequential fashion, \cite{Stolk:An_improved_sweeping_domain_decomposition_preconditioner_for_the_Helmholtz_equation} introduced a modified sweeping pattern, which changes the data dependencies during the sweeps to improve  parallelism.  Finally, \citet{ZepedaDemanet:the_method_of_polarized_traces} introduced the method of polarized traces, which reduces the solver's run-time by leveraging parallelism and fast summation methods.  This paper builds on top of the general framework of the method of polarized traces, which we elaborate on in the sequel.

To date, most studies focus on minimizing the parallel runtime or complexity of a single solve with a single right-hand side.  However, in the scope of seismic inversion, where there can be many thousands of right-hand sides, it is important to consider the overall runtime or complexity of solving \textit{all} right-hand sides.  In this context, linear complexity is $\cO(RN)$, where $N$ is the total number of degrees of freedom (we assume that $N = n^3$ and let $n$ be the number of degrees of freedom in a single dimension of a 3D volume) and $R$ is the number of right-hand sides.

In this paper, we present a solver for the 3D high-frequency Helmholtz equation with a {\it sublinear} online parallel runtime,  given by
\[
\cO( \alpha_{\text{pml}}^2 \max(1,R/L) N \log{N}),
\]
where $N = n^3$ is the total number of unknowns, $L \sim n$ is the number of subdomains in a layered domain decomposition, and $\alpha_{\text{pml}}$ is the number of points needed to implement a high-quality absorbing boundary condition between layers.
We achieve this complexity by comprehensive parallelization of all aspects of the algorithm, including exploiting parallelism in local solves and by pipelining the right-hand sides.
Thus, as long as $R \sim n^2$ (3D), there is a mild $R/L \sim n$ factor impacting the asymptotic complexity.  The solver in this paper is based on the method of polarized traces \citep{ZepedaDemanet:the_method_of_polarized_traces}, a layered domain decomposition method which exploits:
\begin{itemize}
\item  local solvers, using efficient sparse direct solvers at each subdomain,
\item  high-quality transmission conditions between subdomains, implemented via perfectly-matched layers (PML; \cite{Berenger:PML,Johnson:PML}), and
\item  an efficient preconditioner based on polarizing conditions imposed via incomplete Green's integrals.
\end{itemize}
These concepts combine to yield a global iterative method that converges in a small number of iterations.  The method has two stages: an offline stage, that can be precomputed independently of the right-hand sides, and an online stage, that is computed for each right-hand side or by batch processing.

One advantage of the method of polarized traces is that only the degrees of freedom at the interfaces between layers are needed for the bulk of the computation, because the volume problem is reduced to an equivalent surface integral equation (SIE) at the interfaces between layers.  Due to this efficiency, the algorithm requires a smaller memory footprint, which helps make feasible pipelining the right-hand sides.  Pipelining for domain decomposition methods has been previously considered \citep{Stolk:An_improved_sweeping_domain_decomposition_preconditioner_for_the_Helmholtz_equation}, albeit without complexity claims and without a fully tuned communication strategy between subdomains. Moreover, we are unaware of any complexity claims within the context of inversion algorithms, in particular, with focus on full waveform inversion \citep{Taratola:Inversion_of_seismic_reflection_data_in_the_acoustic_approximation}, where there are many right-hand sides that are strongly frequency dependant.

\subsection{Complexity claims}

Suppose that $L$, the number of layers in the domain decomposition, scales as $L \sim n$, i.e., each layer has a constant thickness in number of grid points\footnote{This hypothesis is critical for obtaining a quasi-linear complexity algorithm and is related to the complexity of solving a quasi-2D problem using multi-frontal methods (for further details see \cite{EngquistYing:Sweeping_PML,Poulson_Engquist:a_parallel_sweeping_preconditioner_for_heteregeneous_3d_helmholtz}).}. Each layer is further extended by $\alpha_{\text{pml}}$ grid points in order to implement the PML. It has been documented in \cite{Poulson_Engquist:a_parallel_sweeping_preconditioner_for_heteregeneous_3d_helmholtz} that $\alpha_{\text{pml}}$ needs to grow with problem frequency, $\alpha_{\text{pml}} \sim \log\omega$, in order to obtain a number of Krylov solver iterations to convergence that scales as $\log \omega$.
Additionally, as above, for the 3D problem it is typical that the number of sources $R \sim n^2$, because as frequency and resolution increase, the number of sources in both the in-line and cross-line direction must also increase \citep{Brossier:2D_and_3D_frequency-domain_elastic_wave_modeling_in_complex_media_with_a_parallel_iterative_solver}.

Finally, given that we solve the 3D problem in a high-performance computing (HPC) environment, we assume that the number of computing nodes in the HPC cluster is $\cO(n^3 \log(n)/M)$, with $L \sim n$ layers, $\cO(n^2\log(n)/M)$ nodes inside each layer, and $M$ is the memory of a single node.
The assumptions on node growth come from the fact that computing nodes have finite memory, and thus more nodes are needed to solve larger problems.   As numerical examples will show, using more nodes per layer reduces the runtime per Krylov iteration by enhancing the parallelism of the solves at each subdomain, provided that a carefully designed communication pattern is used, to keep the communication overhead low.

We summarize the asymptotic runtimes in Table \ref{table:complexity}. Like other related methods, the offline stage is linear in $N$ and is independent of the number of right-hand sides. The online runtime is sub-linear, in the sense that linear complexity would be $\cO(RN)$.

\begin{table}
    \begin{center}
        \begin{tabular}{|c|c|}
            \hline
            Stage    & Polarized traces   \\
            \hline
            offline  & $\cO\left(\alpha_{\text{pml}}^3  N \right)$ \\
            \hline
            online   & $\cO(  \alpha_{\text{pml}}^2 \max(1,R/L) N \log{N}) $ \\
            \hline
        \end{tabular}
    \end{center}
    \caption{Runtime of both stages of the algorithm. Note that typically $ \alpha_{\text{pml}} \sim \log \omega$.}\label{table:complexity}
\end{table}

\subsection{Related work} \label{section:related_work}

Modern linear algebra techniques, in particular nested dissection methods \citep{GeorgeNested_dissection} coupled with $\cH$-matrices \citep{Hackbusch:Hierarchical_matrices} have been applied to the Helmholtz problem, yielding, for example: the hierarchical Poincar\'e-Steklov solver \citep{Gillman_Barnett_Martinsson:A_spectrally_accurate_solution_technique_for_frequency_domain_scattering_problems_with_variable_media}, solvers using hierarchical semi-separable (HSS) matrices \citep{Wang:H_multifrontal,Wang_de_Hoop:Massively_parallel_structured_multifrontal_solver_for_time-harmonic_elastic_waves_in_3D_anisotropic_media,Wang_Li_Sia_Situ_Hoop:Efficient_Scalable_Algorithms_for_Solving_Dense_Linear_Systems_with_Hierarchically_Semiseparable_Structures}, or block low-rank (BLR) matrices \citep{Amestoy_Weisbecker:compressed_MUMPS,Amestoy:Fast_3D_frequency_domain_full_waveform_inversion_with_a_parallel_block_low-rank_multifrontal_direct_solver_Application_to_OBC_data_from_the_North_Sea}.

Multigrid methods, once thought to be inefficient for the Helmholtz problem, have been successfully applied to the Helmholtz problem by \citet{Calandra_Grattonn:an_improved_two_grid_preconditioner_for_the_solution_of_3d_Helmholtz} and \citet{Stolk}.  Although these algorithms do not result in a lower computational complexity, their empirical run-times are impressive due to their highly parallelizable nature.  Due to the possibility for efficient parallelization, there has been a renewed interest on multilevel preconditioners such as the one in \cite{Hu_Zhang:Substructuring_Preconditioners_for_the_Systems_Arising_from_Plane_Wave_Discretization_of_Helmholtz_Equations}.

Within the geophysical community, the analytic incomplete LU (AILU) method was explored by \citet{Plessix_Mulder:Separation_of_variable_preconditioner_for_iterativa_Helmholtz_solver} and applied in the context of 3D seismic imaging, resulting in some large computations \citep{Plessix:A_Helmholtz_iterative_solver_for_3D_seismic_imaging_problems}.
A variant of Kazmarc preconditioners \citep{Gordon:A_robust_and_efficient_parallel_solver_for_linear_systems} have been studied and applied to time-harmonic wave equations by \citet{Brossier:2D_and_3D_frequency-domain_elastic_wave_modeling_in_complex_media_with_a_parallel_iterative_solver}.
Although these solvers have, in general, relatively low memory consumption they tend to require many iterations to converge, thus hindering practical run-times.

Domain decomposition methods for solving PDEs date back to \cite{Schwarz:Uber_einen_Grenzubergang_durch_alternierendes_Verfahren}, in which the Laplace equation is solved iteratively (for a more recent treatise, see \citet{Lions:on_the_Schwarz_alternating_method_I}).
The application of domain decomposition to the Helmholtz problem was first proposed by \citet{Despres:domain_decomposition_hemholtz}.
\cite{Cessenat_Despres:application_of_an_ultra_weak_variational_formulation_of_elliptic_pdes_to_the_2_d_helmholtz_problem} further refined this approach with the development of the ultra-weak variational formulation (UWVF) for the Helmholtz equation, in which the basis functions in each element, or sub-domain, are solutions to the local homogeneous equation.
The UWVF approach motivated a series of related methods, such as the partition of unity method of \citet{babuska_melenk:partition_of_unity_method}, the least squared method of \citet{Monk_Wang:A_least-squares_method_for_the_Helmholtz_equation:}, the discontinuous enrichment method by \citet{Farhat:The_discontinuous_enrichment_method}, and Trefftz methods by \citet{Gittelson_Hipmair_Perugia:Trefftz}, and \citet{Perugia:trefft}, among many others. A recent and thorough review of Trefftz and related methods can be found in (\cite{Hiptmair_Moiola_Perugia:A_Survey_of_Trefftz_Methods_for_the_Helmholtz_Equation}).

The results in \citet{Lions:on_the_Schwarz_alternating_method_I} and \citet{Despres:domain_decomposition_hemholtz} have inspired the development of various domain decomposition algorithms, which are now classified as Schwarz algorithms\footnote{For a review on classical Schwarz methods see \citep{Chan:Domain_decomposition_algorithms,Toselli:Domain_Decomposition_Methods_Algorithms_and_Theory}; and for other applications of domain decomposition methods for the Helmholtz equations, see \citep{Bourdonnaye_Farhat_Roux:A_NonOverlapping_Domain_Decomposition_Method_for_the_Exterior_Helmholtz_Problem,Ghanemi98adomain,McInnes_Keyes:Additive_Schwarz_Methods_with_Nonreflecting_Boundary_Conditions_for_the_Parallel_Computation_of_Helmholtz_Problems,Collino:Domain_decomposition_method_for_harmonic_wave_propagation_a_general_presentation,Magoules:Application_of_a_domain_decomposition_with_Lagrange_multipliers_to_acoustic_problems_arising_from_the_automotive_industry,Boubendir:An_analysis_of_the_BEM_FEM_non_overlapping_domain_decomposition_method_for_a_scattering_problem,Astaneh_Guddati:A_two_level_domain_decomposition_method_with_accurate_interface_conditions_for_the_Helmholtz_problem}.}.
However, the convergence rate of such algorithms is strongly dependent on the boundary conditions prescribed at the interfaces between subdomains.
\cite{Gander_Nataf:Optimized_Schwarz_Methods_without_Overlap_for_the_Helmholtz_Equation} introduces an optimal, non-local boundary condition for domain interfaces, which is then approximated by an optimized Robin boundary condition. This last work lead to the introduction of the framework of optimized Schwarz methods in \cite{Gander:Optimized_Schwarz_Methods} to described optimized boundary conditions that provides high convergence.  The design of better interface approximations has been studied in \cite{Gander_Kwok:optimal_interface_conditiones_for_an_arbitrary_decomposition_into_subdomains,Geuzaine:A_quasi-optimal_nonoverlapping_domain_decomposition_algorithm_for_the_Helmholtz_equation,Gander_Zhang:Domain_Decomposition_Methods_for_the_Helmholtz_Equation:_A_Numerical_Investigation,Gander:Optimized_Schwarz_Methods_with_Overlap_for_the_Helmholtz_Equation,Gander:Optimized_Schwarz_Method_with_Two_Sided_Transmission_Conditions_in_an_Unsymmetric_Domain_Decomposition} among many others.

\cite{Engquist_Zhao:Absorbing_boundary_conditions_for_domain_decomposition} introduced absorbing boundary conditions for domain decomposition schemes for elliptic problems and the first application of such techniques to the Helmholtz problem traces back to the AILU factorization (\cite{GanderNataf:ailu_for_hemholtz_problems_a_new_preconditioner_based_on_an_analytic_factorization}).
The sweeping preconditioner, introduced in \cite{EngquistYing:Sweeping_H,EngquistYing:Sweeping_PML}, was the first algorithm to show that those ideas could yield algorithms with quasi-linear complexity. There exists two variants of the sweeping preconditioner which involved using either $\cH$-matrices \citep{EngquistYing:Sweeping_H} or multi-frontal solvers \citep{EngquistYing:Sweeping_PML} to solve the local problem in each thin layer.
These schemes are extended  to different discretizations and physics by \cite{Tsuji_engquist_Ying:A_sweeping_preconditioner_for_time-harmonic_Maxwells_equations_with_finite_elements,Tsuji_Poulson:sweeping_preconditioners_for_elastic_wave_propagation,Tsuji_Ying:A_sweeping_preconditioner_for_Yees_finite_difference_approximation_of_time-harmonic_Maxwells_equations}. Since the introduction of the sweeping preconditioner, several related algorithms with similar claims have been proposed, such as the source transfer preconditioner (\cite{Chen_Xiang:a_source_transfer_ddm_for_helmholtz_equations_in_unbounded_domain,Cheng_Xiang:A_Source_Transfer_Domain_Decomposition_Method_For_Helmholtz_Equations_in_Unbounded_Domain_Part_II_Extensions}), the rapidly converging domain decomposition (\cite{CStolk_rapidily_converging_domain_decomposition}) and its extensions (\cite{Stolk:An_improved_sweeping_domain_decomposition_preconditioner_for_the_Helmholtz_equation}), the double sweep preconditioner (\cite{GeuzaineVion:double_sweep}) and the method of polarized traces (\cite{ZepedaDemanet:the_method_of_polarized_traces}).

\subsection{Organization}
The remainder of this paper is organized as follows: we provide the numerical formulation of the Helmholtz equation, present the reduction to a surface integral equation, and introduce the method of polarized traces for solving the SIE.   Next, we elaborate on the parallelization and communication patterns and examine the empirical complexities and runtimes. Finally, we provide results from several experiments to support our claims.

\section{Formulation} \label{section:formulation}
For this study, we discretize Eq.~\ref{eq:Helmholtz} using the standard second order finite difference method on a regular mesh of $\Omega$, with a grid of size $n_x \times n_y \times n_z$ and a grid spacing $h$.
Note, the method of polarized traces is not restricted to second-order finite differences, but using higher-order finite difference schemes makes the numerical implementation slightly more complicated\footnote{Nor is the method applicable only in the context of finite differences: finite elements \citep{Taus_Demanet_Zepeda:HDG_Helmholtz,Fang_Qian_Zepeda_Zhao:Learning_Dominant_Wave_Directions_For_Plane_Wave_Methods_For_High_Frequency_Helmholtz_Equations} and integral equations \citep{ZepedaZhao:Fast_Lippmann_Schwinger_solver} approaches are also valid.}.
Absorbing  boundary conditions are imposed via perfectly matched layers (PMLs) as described by \citet{Berenger:PML} and \citet{Johnson:PML}.

We describe our PML implementation in detail because the quality and structure of the PML implementation strongly impact the convergence properties of the method.  Following \cite{Berenger:PML,Johnson:PML}, the PML's are implemented via a complex coordinate stretching.
First, we define an extended domain $\hat\Omega$ such that $\Omega \subset \hat\Omega$ and we extend the Helmholtz operator from Eq.~\ref{eq:Helmholtz} to that domain as follows:
\begin{equation} \label{eq:extended_operator_PML}
    \cH = \hat\triangle + m\omega^2 \qquad  \text{ in } \hat\Omega,
\end{equation}
where $m$ is an extension of the squared slowness to $\hat\Omega$ and the extended Laplacian $\hat\triangle$ is constructed by replacing the partial derivatives in the standard Laplacian $\triangle = \partial_{xx} + \partial_{yy} + \partial_{zz}$ with coordinate-stretched partial derivatives defined on $\hat\Omega$:
\begin{equation}
    \partial_x  \rightarrow \beta_x(\x) \partial_x, \,\, \partial_y  \rightarrow \beta_y(\x) \partial_y, \,\, \partial_z  \rightarrow  \beta_z(\x) \partial_z.
\end{equation}
The complex dilation function $\beta_x(\x)$ (and similarly $\beta_y(\x)$ and $\beta_z(\x)$)  is defined as
\begin{equation} \label{appendix:eq:def_alpha}
    \beta_x(\x) =  \frac{1}{1+ i \frac{\sigma_x(\x)}{ \omega } },
\end{equation}
where the PML profile function $\sigma_x(\x)$ (and similarly $\sigma_y(\x)$ and $\sigma_z(\x)$) is,
\begin{equation}
    \sigma_x(\x) =  \left \{\begin{array}{rl}
                        \frac{C}{\delta_{\text{pml}}} \left (\frac{x } {\delta_{\text{pml}}} \right)^2,         & \text{if } x \in (-\delta_{\text{pml}}, 0 ),\\
                        0               ,                                                                           & \text{if } x \in [ 0, L_x ],  \\
                        \frac{C}{\delta_{\text{pml}}} \left (\frac{x - L_x }{\delta_{\text{pml}}}\right)^2,     & \text{if } x \in ( L_x , L_x + \delta_{\text{pml}} ),\\
                            \end{array} \right .
\end{equation}
where $L_x$ is the length of $\Omega$ in the $x$ dimension and $\delta_{\text{pml}}$ is the length of the extension.
In general, $\delta_{\text{pml}}$ grows slowly with the frequency, i.e., $\delta_{\text{pml}} \propto \cO(\log{\omega})$, in order to obtain enough absorption as the frequency increases.
The constant $C$ is chosen to provide enough absorption.  In practice, $\delta_{\text{pml}}$ and $C$ can be seen as parameters to be tuned for accuracy versus efficiency.

The extended Helmholtz operator provides the definition of the global continuous problem,
\begin{equation} \label{eq:Helmholtz_pml}
\cH u = f_s, \qquad \text{in } \hat\Omega,
\end{equation}
which is then discretized using finite differences to obtain the discrete global problem,
\begin{equation}
    \mathbf{H} \u = \f_s. \label{eq:discrete_Helmholtz}
\end{equation}

In the method of polarized traces, $\Omega$ is decomposed into a set of $L$ layers, $\{\Omega^{\ell}\}_{\ell=1}^{L}$.  Without loss of generality, we assume that the decomposition is in the $z$ dimension.  Each subdomain $\Omega^{\ell}$ is extended to include an absorbing region, as above, yielding the extended subdomain $\hat\Omega^{\ell}$.  For boundaries of $\hat\Omega^{\ell}$ shared with $\hat\Omega$, the absorbing layer is considered to be inherited from the global problem.  For the intra-layer boundaries of $\Omega^{\ell}$, i.e., those due to the partitioning of $\Omega$, the extension to the additional absorbing layers in $\hat\Omega^{\ell}$ are necessary to prevent reflections at layer interfaces which are detrimental to convergence.

\begin{figure}[h]
    \begin{center}
        \includegraphics[trim = 0mm 0mm 0mm 0mm, clip, width=8.5cm]{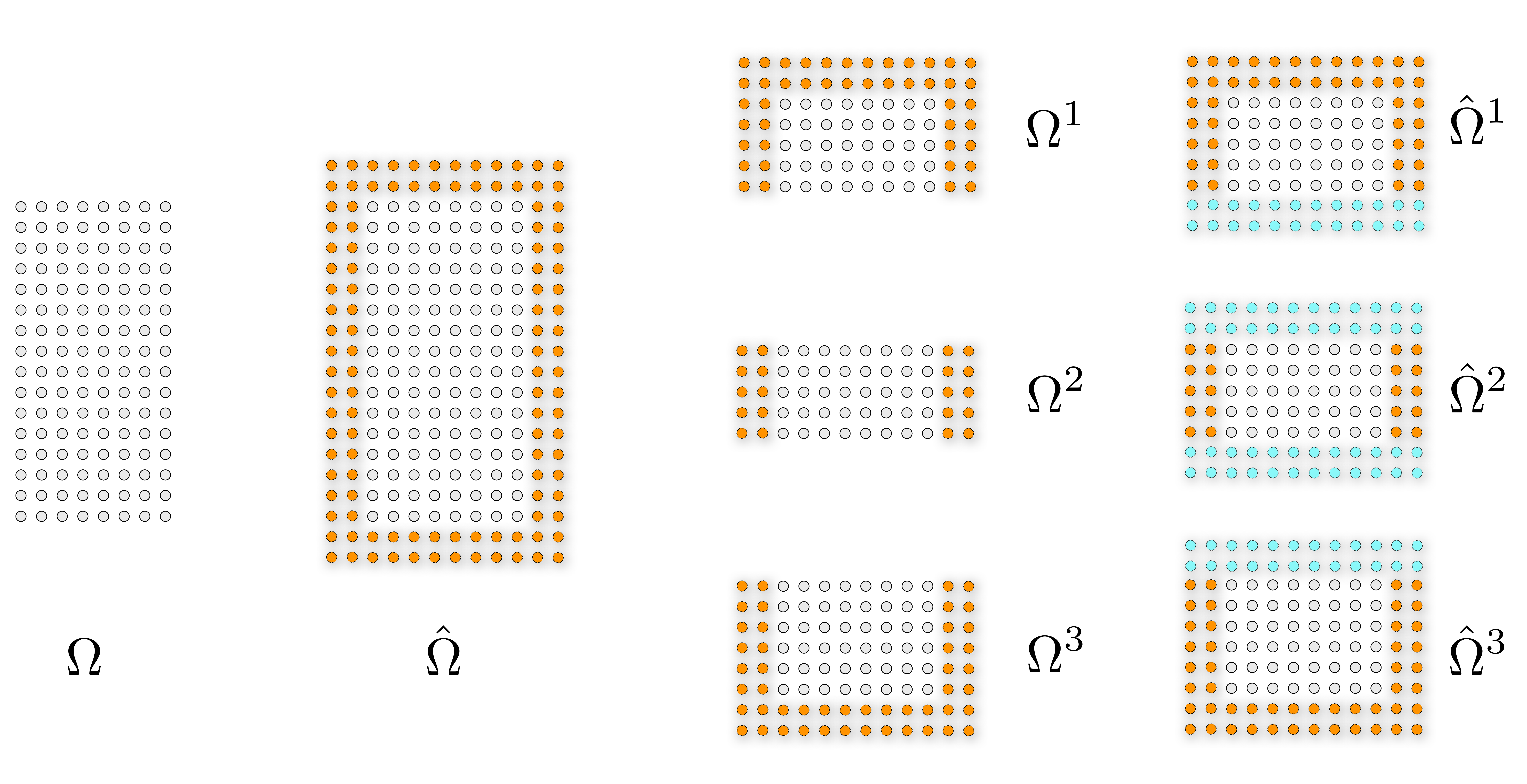}
    \end{center}
    \caption{Sketch in 2D of the partition of the domain in layers. The domain $\Omega$ is extended to $\hat\Omega$ by adding the PML nodes (orange). After decomposition into subdomains $\Omega^{\ell}$, the internal boundaries are padded with extra PML nodes (light blue) resulting in the subdomains $\hat\Omega^{\ell}$.} \label{fig:DDM_sketch_PML}
\end{figure}

Then, the local Helmholtz problem is
\begin{equation}
   \cH^{\ell} v^{\ell} := \triangle^{\ell} v^{\ell} (\x)  + m^{\ell} \omega^2 v^{\ell} (\x)= f_s^{\ell}(\x) \qquad \text{in } \hat\Omega^{\ell}, \label{eq:local_Helmholtz}
\end{equation}
where $m^{\ell}$ and $f^{\ell}$ are the local restrictions of the model parameters and source functions to $\hat\Omega^{\ell}$, generated by extending $m \chi_{\Omega^{\ell}}$ and $f \chi_{\Omega^{\ell}}$ to $\hat\Omega^{\ell}$, where $\chi_{\Omega^{\ell}}$ is the characteristic function of $\Omega^{\ell}$.   The local Laplacian $\triangle^{\ell}$ is defined using the same coordinate stretching approach as above, except on $\hat\Omega^{\ell}$.  As before, $\delta_{\text{pml}}$, and thus $\alpha_{\text{pml}}$, must scale as $\log{\omega}$ to obtain the convergence rate claimed in this paper.

We discretize the local problem in Eq.~\ref{eq:local_Helmholtz} resulting in the discrete local Helmholtz system
\begin{equation}
    \mathbf{H}^{\ell} \v^{\ell} = \f_s^{\ell}. \label{eq:discrete_local_Helmholtz}
\end{equation}
For the finite difference implementation in this paper, we assume a structured, equispaced Cartesian mesh with mesh points $\x_{i,j,k} = (x_i,y_j, z_k) = (ih,jh, kh)$.
Assuming the same ordering \citet{ZepedaDemanet:the_method_of_polarized_traces}, we write the global solution in terms of the depth index,
\begin{equation} \label{eq:global_ordering}
\u = (\u_1, \u_2, ..., \u_{n_z}),
\end{equation}
where $\u_k$ is a plane sampled at constant depth $z_k$, or in MATLAB notation,
\begin{equation} \label{eq:trace_ordering}
\u_k = (u_{:,:,k}).
\end{equation}

Let $\u^{\ell}$ be the local restriction of $\u$ to $\Omega^{\ell}$, i.e., $\u^{\ell}  = \chi_{\Omega^{\ell}}\u$.
Following the above notation, $\u^{\ell}_k$ is the local solution trace in the plane at local depth $z_k^{\ell}$.
For notational convenience, we renumber the local depth indices so that $\u^{\ell}_1$ and $\u^{\ell}_{n^{\ell}}$ are the top and bottom planes of the bulk domain.  Points due to the PML are not considered\footnote{With this renumbering the local depth index $z_k^{\ell}$ maps to the global depth index $z_{n_c^{\ell}+k}$ where $n_c^{\ell} = \sum_{j=1}^{\ell-1} n^j$.}.
Finally, let
\begin{align} \label{eq:traces_definition}
    \underline{\u} = \left (\u^{1}_{n^1} , \u^{2}_{1}, \u^{2}_{n^2}, ...,\u^{L-1}_{1}, \u^{L-1}_{n^{L-1}}, \u^{L}_{1}  \right)^{t},
\end{align}
be the vector of interface traces for all $L$ layers.

To map solution vectors at fixed depth planes back to the discretized whole volume of $\Omega_{\ell}$, we define the Dirac delta at a fixed depth,
\begin{equation}
    (\delta(z-z_p) \v_q)_{i,j,k} = \left \{   \begin{array}{cl}
                                                0                                 & \text{if }  k \neq p, \\
                                                \frac{(\v_q)_{i,j}}{h^3} & \text{if }  k = p.
                                            \end{array}
                                 \right .
\end{equation}
This definition of the numerical Dirac delta is specific to a classical finite difference discretization.  If the discretization changes, it is still possible to define a numerical Dirac delta  using the approach developed in \cite{ZepedaDemanet:Nested_domain_decomposition_with_polarized_traces_for_the_2D_Helmholtz_equation}.

\subsection{Accuracy}

Solving the Helmholtz equation in the high-frequency regime is notoriously because:
\begin{itemize}
\item it is difficult to efficiently discretize the PDE, and
\item the resulting linear system is difficult to solve in a scalable and efficient fashion.
\end{itemize}
In this paper, we focus on the second issue.  However, for completeness, we provide a brief overview of difficulties associated to the discretization.

From the Shannon-Nyquist sampling theorem, an oscillatory function at frequency $\omega$ requires $\cO(\omega^d)$ degrees of freedom to be accurately represented, without aliasing. For example, to accurately represent the solution of Eq.~\ref{eq:Helmholtz}, only $\cO(\omega^3)$ degrees of freedom are required.   Obviously, accuracy is still limited by the error in the discretization of the differential operator.
Even if the medium is very smooth, standards methods based on finite differences and finite elements are subject to pollution error, i.e., the ratio between the error of the numerical approximation and the best approximation cannot be bounded independently of $\omega$ \citep{Babuska:A_Generalized_Finite_Element_Method_for_solving_the_Helmholtz_equation_in_two_dimensions_with_minimal_pollution,Ihlenburg_Babuska:Finite_element_solution_of_the_Helmholtz_equation_with_high_wave_number_Part_I:_The_h_version_of_the_FEM,Babuska:Stable_Generalized_Finite_Element_Method}.

The direct consequence of pollution error is that the approximation error, i.e., the error between the analytical and the numerical solution to the linear system, increases with the frequency, even if $n \sim \omega$. Thus, to obtain a bounded approximation error independent of the frequency, it is required to oversample the wavefield, relative to the Shannon-Nyquist criterion, i.e., $n$ needs to grow faster than $\omega$. Unfortunately, oversampling provides discretizations with a suboptimal number of degrees of freedom with respect to the frequency.
To alleviate pollution error, several new approaches have been proposed, which can be broadly classified into two groups:
\begin{enumerate}
\item methods using standard polynomial bases with modified variational formulations \citep{GOLDSTEIN:The_weak_element_method_applied_to_Helmholtz_type_equations,Melenk_Sauter:Wavenumber_Explicit_Convergence_Analysis_for_Galerkin_Discretizations_of_the_Helmholtz_Equation,Melenk_Parsania_Sauter_13:General_DG_Methods_for_Highly_Indefinite_Helmholtz_Problems,Moiola_Spence:Is_the_Helmholtz_Equation_Really_Sign_Indefinite,Graham_Lohndorf_Melenk_Spence:When_is_the_error_in_the_h-BEM_for_solving_the_Helmholtz_equation_bounded_independently_of_k};

\item methods based on well known variational formulations but using non-standard basis, such as plane waves \citep{Hiptmair_Moiola_Perugia:A_Survey_of_Trefftz_Methods_for_the_Helmholtz_Equation,Perugia_Pietra_Russo:PW_VEM}, polynomials modulated by plane waves \citep{Betcke_Phillips:Approximation_by_dominant_wave_directions_in_plane_wave_methods,DBLP:journals/jcphy/NguyenPRC15}, or other specially adapted functions.

\end{enumerate}

Even though the methods mentioned above have been successful in reducing pollution error, the resulting linear systems cannot, in general, be solved in quasi-linear time or better because the matrices either have a high degree of interconnectivity or are extremely ill-conditioned.
However, some new fast algorithms have recently been proposed for solving the Helmholtz equation without pollution error with quasi-linear complexity for media that are homogeneous up to smooth and compactly supported heterogeneities  (\citet{ZepedaZhao:Fast_Lippmann_Schwinger_solver,Fang_Qian_Zepeda_Zhao:Learning_Dominant_Wave_Directions_For_Plane_Wave_Methods_For_High_Frequency_Helmholtz_Equations}; and references therein).
In the case of highly heterogeneous media the accuracy of finite elements has not been extensively studied, although methods of efficient discretizations for highly heterogeneous media, coupled with fast algorithms are emerging \citep{Taus_Demanet_Zepeda:HDG_Helmholtz}.

In this paper, we assume that waves will propagate in very general and highly heterogeneous media, thus, we do not have a theoretical framework to assess the accuracy. Instead, we use numerical experiments to check the accuracy of the solution. Our numerical experiments show, for the cases considered in this paper, using 10 points per wavelength results in roughly 1 digit of accuracy at the highest frequency considered.

\subsection{Reduction to a surface integral equation}

The global solution is related to the local layer solutions by coupling the subdomains using the Green's representation formula (GRF) within each layer.
The resulting surface integral equation (SIE), posed at the interface between layers, effectively reduces the problem from the global domain $\Omega$ to the interfaces between layers. The resulting SIE has the form
\begin{equation}
    \underline{\mathbf{M}} \underline{\u} = \underline{\f}, \label{eq:SIE}
\end{equation}
where $\underline{\mathbf{M}}$ is formed by interface-to-interface Green's functions, $\underline{\u}$ is defined in Eq.~\ref{eq:traces_definition}, and $\underline{\f}$ is the right-hand-side, formed as in Line 8 of Alg.~\ref{alg:SIE_solver}.

\begin{figure}[t]
    \begin{center}
        \includegraphics[trim = 20mm 22mm 16mm 17mm, clip, width=4cm]{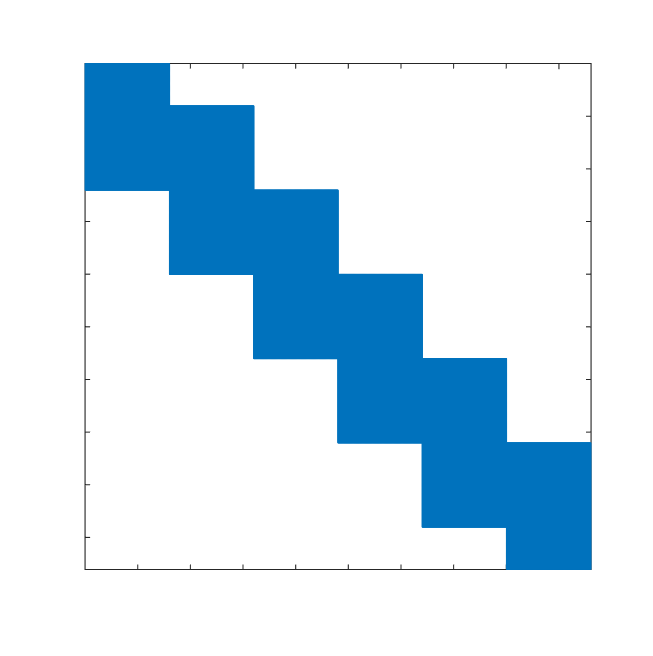}  \includegraphics[trim = 20mm 22mm 16mm 17mm, clip, width=4cm]{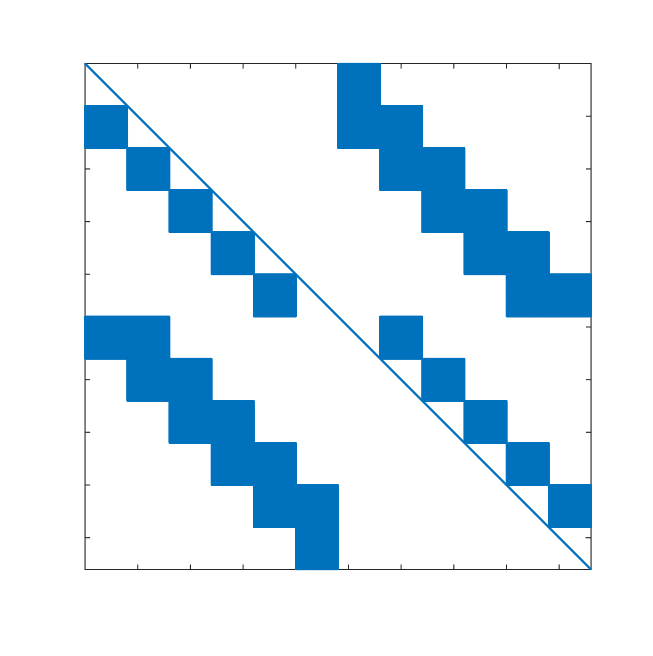}
    \end{center}
    \caption{Sparsity pattern of the SIE matrix in Eq.~\ref{eq:SIE} (left), and the polarized SIE matrix in Eq.~\ref{eq:polarized_SIE} (right) .} \label{fig:M_polarized}
\end{figure}

The matrix $\underline{\mathbf{M}}$ is a block banded matrix (Fig.~\ref{fig:M_polarized}, left) of size $2(L-1)n^2 \times 2(L-1)n^2$. Theorem 1 of \citet{ZepedaDemanet:the_method_of_polarized_traces} gives that the solution of Eq.~\ref{eq:SIE} is exactly the restriction of the solution of Eq.~\ref{eq:discrete_Helmholtz} to the interfaces between layers.

Following \cite{ZepedaDemanet:the_method_of_polarized_traces}, if the traces of the exact solution are known, then it is possible to apply the GRF to locally reconstruct exactly the global solution within each layer.
Equivalently, the reconstruction can be performed by modifying the local source with a measure supported on the layer interfaces and solving the local system with the local solver, as seen in lines 11-12 of Alg.~\ref{alg:SIE_solver},  where a high-level sketch of the algorithm to solve the 3D high-frequency Helmholtz equation is given.

To efficiently solve the 3D problem, it is critical that the matrix $\underline{\mathbf{M}}$ is never explicitly formed.
Instead, a matrix-free approach \citep{ZepedaDemanet:Nested_domain_decomposition_with_polarized_traces_for_the_2D_Helmholtz_equation} is used to apply the blocks of $\underline{\mathbf{M}}$ via applications of the local solver, using equivalent sources supported at the interfaces between layers, as shown in Alg.~\ref{alg:applyM}. Moreover, as seen in Alg.~\ref{alg:applyM}, the application of $\underline{\mathbf{M}}$ is easily implemented in parallel, with a small communication overhead.
The only non-embarrassingly parallel stage of Alg.~\ref{alg:SIE_solver} is the solution of Eq.~\ref{eq:SIE}, which is inherently sequential.

Given that $\underline{\mathbf{M}}$ is never explicitly formed, an iterative method is the natural choice for solving Eq.~\ref{eq:SIE}.
In practice, the condition number of $\underline{\mathbf{M}}$  is very large and it has a wide spectrum is the complex plane, which implies that a large number of iterations are required to achieve convergence.  To alleviate this problem, we apply the method of polarized traces, as an efficient preconditioner for Eq.~\ref{eq:SIE}, which we describe below.

\begin{algorithm} Offline computation \label{alg:offline_computation}
    \begin{algorithmic}[1]
        \Function{ $[  L^{\ell}, U^{\ell} ]$ = Factorization}{$m$, $\omega$}
            \For{  $\ell = 1: L$ }
                \State$\mathbf{H}^{\ell} = \triangle^{\ell}  + m^{\ell} \omega^2 $ \Comment{Build the local system}
                \State $\mathbf{L}^{\ell} \mathbf{U}^{\ell} = \mathbf{H}^{\ell} $  \Comment{Compute the LU factorization}
            \EndFor
        \EndFunction
    \end{algorithmic}
  \end{algorithm}

\begin{algorithm} Online computation using the SIE reduction \label{alg:SIE_solver}
    \begin{algorithmic}[1]
        \Function{ $\u$ = Helmholtz solver}{ $\mathbf{f}$ }
            \For{  $\ell = 1: L$ }
                \State $ \mathbf{f}^{\ell} = \mathbf{f}\chi_{\Omega^{\ell}} $                           \Comment{partition the source}
            \EndFor
            \For{  $\ell = 1: L$ }
                \State $\v^{\ell} = (\mathbf{H}^{\ell})^{-1} \mathbf{f}^{\ell}  $                \Comment{solve local problems}
            \EndFor
            \State $ \underline{\mathbf{f}} =  \left ( \v^1_{n^1}, \v^2_1 ,\v^2_{n^2} ,\hdots ,\v^L_1 \right )^{t} $           \Comment{form r.h.s.}
            \State $\underline{\u} = \left( \underline{\mathbf{M}} \right )^{-1} \underline{\f} $                                                             \Comment{solve Eq.~\ref{eq:SIE}}
            \For{  $\ell = 1: L$ }
                \State $ \begin{array}{ll} \mathbf{g}^{\ell} =  & \f^{\ell} + \delta(z_{1}-z)\u^{\ell-1}_{n^{\ell-1}} - \delta(z_{0}-z)\u^{\ell}_{1}  \\
                                                                & - \delta(z_{n^{\ell}+1}-z)\u^{\ell}_{n^{\ell}}  + \delta(z_{n^{\ell}}-z)\u^{\ell+1}_{1}
                         \end{array}$
                \State $ \u^{\ell}  = (\mathbf{H}^{\ell})^{-1} \mathbf{g}^{\ell}$ \Comment{inner solve}
            \EndFor
            \State $\u = (\u^1 , \u^2, \hdots, \u^{L-1}, \u^{L})^t $                                                                                                                \Comment{concatenate}
        \EndFunction
    \end{algorithmic}
  \end{algorithm}

\begin{algorithm} Application of the boundary integral matrix $\underline{\mathbf{M}}$ \label{alg:applyM}
    \begin{algorithmic}[1]
        \Function{ $\underline{\u}$ = Boundary Integral}{ $\underline{\v}$ }
            \State $ \tilde{\f}^{1} =  - \delta(z_{n^{1}+1}-z)\v^{1}_{n^{\ell}}  + \delta(z_{n^{1}}-z)\v^{2}_{n^{1}}  $
            \State $ \mathbf{w}^{1}  = (\mathbf{H}^{1})^{-1} \tilde{\f}^{1}$
            \State $\u^{\ell}_{n^{\ell}}   = \mathbf{w}^{\ell}_{n^{\ell}} - \v^{\ell}_{n^{\ell}} $
            \For{  $\ell = 2: L-1$ }
                \State $ \begin{array}{ll} \tilde{\f}^{\ell} =  &   \delta(z_{1}-z)\v^{\ell-1}_{n^{\ell-1}} - \delta(z_{0}-z)\v^{\ell}_{1}  \\
                                                                & - \delta(z_{n^{\ell}+1}-z)\v^{\ell}_{n^{\ell}}  + \delta(z_{n^{\ell}}-z)\v^{\ell+1}_{1}
                         \end{array}$
                \State $ \mathbf{w}^{\ell}  = (\mathbf{H}^{\ell})^{-1} \tilde{\f}^{\ell} $ \Comment{inner solve}
                \State $\u^{\ell}_{1}          = \mathbf{w}^{\ell}_{1}        - \v^{\ell}_{1}; \qquad \u^{\ell}_{n^{\ell}}   = \mathbf{w}^{\ell}_{n^{\ell}} - \v^{\ell}_{n^{\ell}} $
            \EndFor
            \State $ \tilde{\f}^{L} = \delta(z_{1}-z)\v^{L-1}_{n^{L-1}} - \delta(z_{0}-z)\v^{L}_{1}$
            \State $ \mathbf{w}^{L}  = (\mathbf{H}^{L})^{-1} \tilde{\f}^{L}$
            \State $\u^{L}_{1}          = \mathbf{w}^{L}_{1}        - \v^{L}_{1}$
        \EndFunction
    \end{algorithmic}
\end{algorithm}

\section{Method of Polarized traces}

Reducing the Helmholtz problem to a SIE allows us to efficiently parallelize most of the computation required to solve Eq.~\ref{eq:discrete_Helmholtz}.
The only remaining sequential bottleneck is the solution of Eq.~\ref{eq:SIE}.
Given the size and the distributed nature of  $\underline{\mathbf{M}}$, iterative methods, such as GMRES (\cite{Saad_Schultz:GMRES}) or Bi-CGSTAB (\cite{van_der_Vorst:BiCGSTAB}), are the logical approach for solving Eq.~\ref{eq:SIE}.
However, numerical experiments indicate that the condition number of  $\underline{\mathbf{M}}$ scales as $\cO(h^{-2})$, or as $\cO(\omega^2)$ in the high-frequency regime~\citep{ZepedaDemanet:the_method_of_polarized_traces}.
The number of iterations required for schemes like GMRES to converge is proportional to the condition number of the system, yielding poor scalability for solving the SIE at high frequencies.
To alleviate this problem, we use the method of polarized traces to convert the SIE to an equivalent problem, which is easily preconditioned.
This preconditioned system only requires $\cO(\log \omega)$ GMRES iterations\footnote{This scaling is empirically deduced , under the assumption that no large resonant cavities are present in the media.}, i.e., it is comparatively  independent of the frequency.
Here, we provide a high-level review of the method of polarized traces and its implementation and we direct the reader to \citet{ZepedaDemanet:the_method_of_polarized_traces} for a detailed exposition.

\subsection{Preconditioner}

As seen in the previous discussion, Eq.~\ref{eq:SIE} is the result of decomposing the domain into a set of layers and reducing the Helmholtz problem to an equivalent SIE on the interfaces between the subdomains.
To precondition the SIE with the method of polarized traces, the solution at the interfaces is decomposed in up- and down-going components such that
\begin{equation}
   \underline{\u}  = \underline{\u}^{\uparrow} + \underline{\u}^{\downarrow},
\end{equation}
which defines the polarized wavefield
\begin{equation}
   \underline{\underline{\u}}  = \left (\begin{array}{c}
                                            \underline{\u}^{\downarrow} \\
                                            \underline{\u}^{\uparrow}
                                        \end{array}
                                 \right ).
\end{equation}
By introducing the polarized wavefield, we have deliberately doubled the unknowns and produced an underdetermined system.
To close the system, we impose annihilation, or polarizing, conditions (see Section 3 of \cite{ZepedaDemanet:the_method_of_polarized_traces}) that are encoded in matrix form as
\begin{equation}
    \underline{\mathbf{A}}^{\uparrow} \underline{\u}^{\uparrow}=0, \,\,\, \text{ and } \,\,\, \underline{\mathbf{A}}^{\downarrow} \underline{\u}^{\downarrow}=0. \label{eq:annihilation_conditions}
\end{equation}


Requiring that the solution satisfies both Eq.~\ref{eq:SIE} and the annihilation conditions yields another equivalent formulation,
\begin{equation}
    \underline{\underline{\mathbf{M}}} \, \underline{\underline{\u}} = \underline{\underline{\f}}_s, \label{eq:polarized_SIE}
\end{equation}
where
\begin{equation}
    \underline{\underline{\mathbf{M}}}  =   \left [ \begin{array}{cc}
                                                        \underline{\mathbf{M}} & \underline{\mathbf{M}} \\
                                                        \underline{\mathbf{A}}^{\downarrow} & \underline{\mathbf{A}}^{\uparrow}
                                                    \end{array}
                                            \right ], \,\, \text{and} \,\,
    \underline{\underline{\f}}_s  = \left ( \begin{array}{c}
                                                \underline{\mathbf{f}}_s \\
                                                0
                                            \end{array}
                                    \right ).
\end{equation}
Following a series of basic algebraic operations and permutations (see \cite{ZepedaDemanet:the_method_of_polarized_traces} for the full details), we obtain an equivalent formulation of the polarized SIE matrix in Eq.~\ref{eq:polarized_SIE}, given by
\begin{equation}
    \underline{\underline{\mathbf{M}}}  =   \left [ \begin{array}{cc}
                                                        \underline{\mathbf{D}}^{\downarrow} & \underline{\mathbf{U}} \\
                                                        \underline{\mathbf{L}}              & \underline{\mathbf{D}}^{\uparrow}
                                                    \end{array}
                                            \right ].
\end{equation}
There exist straightforward, parallel algorithms for applying the block matrices $\underline{\mathbf{D}}^{\downarrow}$, $\underline{\mathbf{D}}^{\uparrow}$, $\underline{\mathbf{L}}$, and $\underline{\mathbf{U}}$. By construction $\underline{\mathbf{D}}^{\downarrow}$ and $\underline{\mathbf{D}}^{\uparrow}$ can be easily inverted using block forward and backward substitution because they are block triangular with identity blocks on their diagonals.
The blocks that appear in the sparsity pattern of $\underline{\underline{\mathbf{M}}}$ (Fig.~\ref{fig:M_polarized}; right) are a direct manifestation of interactions between the layer interfaces.

While the resulting block linear system can be solved using standard matrix-splitting iterations, such as block Jacobi iteration or block Gauss-Seidel iteration \citep{Saad:iterative_methods_for_sparse_linear_systems}, it is natural to continue to use GMRES to solve the system due to the parallel nature of applying the constituent blocks of $\underline{\underline{\mathbf{M}}}$.
However, the structure of $\underline{\underline{\mathbf{M}}}$ is convenient for using a single iteration of Gauss-Seidel as a preconditioner,
\begin{equation}
    \mathbf{P}\underline{\underline{\mathbf{M}}} \, \underline{\underline{\u}} = \mathbf{P}\underline{\underline{\f}}_s, \label{eq:precond_polarized_SIE}
\end{equation}
where the preconditioning matrix is
\begin{equation}  \label{eq:preconditioner_GS}
\mathbf{P} =  \left(
                                \begin{array}{ll}
                                    \mathbf{\underline D}^\downarrow & \mathbf{O}          \\
                                    \mathbf{\underline L}            & \mathbf{\underline D}^\uparrow \\
                                \end{array}
                            \right)^{-1}.
\end{equation}
In the subsequent sections, we will elaborate on the physical and numerical meanings of the constituent blocks of $\underline{\underline{\mathbf{M}}}$ and $\mathbf{P}$.

\subsubsection{Polarization}

The main novelty of the method of polarized traces is due to the polarization conditions, which are encoded in the matrices $\mathbf{\underline{A}}^{\uparrow}$ and  $\mathbf{\underline{A}}^{\downarrow}$.
The polarizing conditions provide a streamlined way to define an iterative solver using standard matrix splitting techniques, and thus an efficient preconditioner for Krylov methods, such as GMRES.

The polarization conditions are constructed by projecting the solution on two orthogonal sets, physically given by waves traveling upwards and downwards.
Similar constructs are well-known to the geophysics community, as methods that decompose wavefields into distinct down- and up-going components are the backbone of several imaging techniques (see \cite{Zhang:The_Theory_of_True_Amplitude_One-Way_Wave_Equation_Migrations} and references therein).
Commonly, the decomposition is  obtained using discretizations of pseudo-differential operators, which can be interpreted as separating the wavefield into a set of wave-atoms traveling in the different directions which are then propagated accordingly.
Methods for decomposing and locally extrapolating directionally decomposed wavefields are well documented (\cite{Wu:Wide-angle_elastic_wave_one-way_propagation_in_heterogeneous_media_and_an_elastic_wave_complex_screen_method,Collino_Joly:Splitting_of_Operators_Alternate_Directions_and_Paraxial_Approximations_for_the_Three_Dimensional_Wave_Equation,Ristow:3D_implicit_finite-difference_migration_by_multiway_splitting,deHoop:generalization_of_the_phase_screen_approximation_for_the_scattering_of_acoustic_waves}).

In our case, we rewrite the decomposition condition as an integral relation between the Neumann and Dirichlet data of the wavefield, which ultimately leads to the annihilation conditions in Eq.~\ref{eq:annihilation_conditions}.
The pair composed of the Neumann and Dirichlet traces should lie within the null space of an integral operator defined on an interface, which allows the decomposition of the total wavefield into the up- and down-going components, with each having a clear physical interpretation.
In particular, an up-going wavefield is a wavefield generated by a source located beneath the interface and it satisfies a radiation condition at positive infinity and a down-going wavefield is a wavefield generated by a source located above the interface and it satisfies a radiation condition at negative infinity.
As detailed in \cite{ZepedaDemanet:the_method_of_polarized_traces}, defining the decomposition in this manner allows us to extrapolate each component in a stable manner using an incomplete Green's integral.

The extrapolation of up-going components is performed algorithmically by the inversion of $\underline{\mathbf{D}}^{\uparrow}$ and in the same fashion the extrapolation of down-going components is performed by the inversion of $\underline{\mathbf{D}}^{\downarrow}$.
Moreover, the application of the operator $\underline{\mathbf{L}}$ isolates the up-going reflections due to down-going waves interacting with the material in each subdomain, and similarly for the operator $\underline{\mathbf{U}}$.

The application of the preconditioner to a decomposed wavefield,
\begin{equation}
\mathbf{P}   \left(  \begin{array}{c} \underline{\v}^{\downarrow} \\
                                      \underline{\v}^{\uparrow}
                        \end{array}
                \right) =   \left(
                                \begin{array}{c}(\underline{\mathbf{D}}^{\downarrow})^{-1} \underline{\v}^{\downarrow} \\
                                                (\underline{\mathbf{D}}^{\uparrow})^{-1} \mathbf{r}
                                \end{array}
                            \right ),
\end{equation}
for $\mathbf{r} = \left ( \v^{\uparrow} -\underline{\mathbf{L}}  (\underline{\mathbf{D}}^{\downarrow})^{-1} \underline{\v}^{\downarrow} \right)$,
can be physically interpreted as follows:
\begin{enumerate}
    \item  $(\underline{\mathbf{D}}^{\downarrow})^{-1} \underline{\v}^{\downarrow}$: extrapolate the down-going components by propagating them downwards,
    \item  $\mathbf{r} = \left ( \v^{\uparrow} -\underline{\mathbf{L}}  (\underline{\mathbf{D}}^{\downarrow})^{-1} \underline{\v}^{\downarrow} \right)$: compute the local reflection of the extrapolated field and add them to the up-going components,
    \item $(\underline{\mathbf{D}}^{\uparrow})^{-1} \mathbf{r}$: extrapolate the up-going components by propagating them upwards.
\end{enumerate}

\subsection{Algorithms}

As with the application of $\underline{\mathbf{M}}$ in Alg.~\ref{alg:applyM}, we construct matrix-free methods for solving $ ( \underline{\mathbf{D}}^{\downarrow})^{-1}$ and $ (\underline{\mathbf{D}}^{\uparrow})^{-1}$  (Algs.~\ref{alg:downwardsSweep} and \ref{alg:upwardsSweep}), where local solves are applied in an inherently sequential fashion.  To complete the preconditioner, a matrix-free (and embarrassingly parallel) algorithm for applying $\underline{\mathbf{L}}$ is given in Alg~\ref{alg:upwardsReflections}.  Similar algorithms for applying $\underline{\mathbf{U}}$, $\underline{\mathbf{D}}^{\uparrow}$, and $\underline{\mathbf{D}}^{\downarrow}$, as well as the complete matrix-free algorithm for applying $\underline{\underline{\mathbf{M}}}$, are provided in the Appendix.

In solving the systems for $ ( \underline{\mathbf{D}}^{\downarrow})^{-1}$ and $ (\underline{\mathbf{D}}^{\uparrow})^{-1}$, each application of the local solver is local to each layer, which means that some communication is required to transfer solution updates from one layer to the next. The sequential nature of the method for solving these systems implies that only one set of processors, those assigned to the current layer, are working at any given stage of the algorithm.
This is illustrated in Fig.~\ref{fig:nodes_load}, where each block represents a local solve and the execution path moves from left to right.  As explained in \cite{ZepedaDemanet:the_method_of_polarized_traces}, it is possible to apply $\underline{\mathbf{D}^{\downarrow}}$ and $\underline{\mathbf{L}}$ simultaneously, thus decreasing the number of local solves per layer.

\begin{algorithm} Downward sweep, application of $ ( \mathbf{\underline{D}}^{\downarrow} )^{-1}$ \label{alg:downwardsSweep}
    \begin{algorithmic}[1]
        \Function{ $\underline{\u}^{\downarrow}$ = Downward Sweep}{ $\underline{\v}^{\downarrow}$ }
            \State $\u^{\downarrow,1}_{n^1}   = -\v^{\downarrow,1}_{n^1}$
            \State $\u^{\downarrow,1}_{n^1+1} = -\v^{\downarrow,1}_{n^1+1} $
            \For{  $\ell = 2: L-1$ }
                \State $ \tilde{\f}^{\ell} = \delta(z_{1}-z)\u^{\downarrow,\ell-1}_{n^{\ell-1}}  - \delta(z_{0}-z)\u^{\downarrow,\ell-1}_{n^{\ell-1}+1} $
                \State $ \mathbf{w}^{\ell}  = (\mathbf{H}^{\ell})^{-1} \tilde{\f}^{\ell} $
                \State $\u^{\downarrow,\ell}_{n^{\ell}}   = \mathbf{w}_{n^{\ell}} - \v^{\downarrow,\ell}_{n^{\ell}}$
                \State $\u^{\downarrow,\ell}_{n^{\ell}+1} = \mathbf{w}_{n^{\ell}+1}-\v^{\downarrow,\ell}_{n^{\ell}+1} $
            \EndFor
            \State $\underline{\u}^{\downarrow} =  \left (\u^{\downarrow,1}_{n^1} , \u^{\downarrow,1}_{n^1+1}, \u^{\downarrow,2}_{n^2}, ..., \u^{\downarrow,L-1}_{n^{L-1}}, \u^{\downarrow,L-1}_{n^{L-1}+1}  \right)^{t} $
        \EndFunction
    \end{algorithmic}
\end{algorithm}

\begin{algorithm} Upward sweep, application of $ ( \mathbf{\underline{D}}^{\uparrow}  )^{-1}$  \label{alg:upwardsSweep}
    \begin{algorithmic}[1]
        \Function{ $\underline{\u}^{\uparrow}$ = Upward sweep}{ $\underline{\v}^{\uparrow}$ }
            \State $\u^{\uparrow,L}_{0}   = -\v^{\uparrow,L}_{0}$
            \State $\u^{\uparrow,L}_{1} = -\v^{\uparrow,L}_{1} $
            \For{  $\ell = L-1:2$ }
                \State $ \tilde{\f}^{\ell} = - \delta(z_{n^{\ell}+1}-z)\u^{\uparrow,\ell+1}_{0}  + \delta(z_{n^{\ell}}-z)\u^{\uparrow,\ell+1}_{1} $
                \State $ \mathbf{w}^{\ell}  = (\mathbf{H}^{\ell})^{-1} \tilde{\f}^{\ell}  $
                \State $\u^{\uparrow,\ell}_{1}   = \mathbf{w}^{\ell}_{1} - \v^{\uparrow,\ell}_{1}$
                \State $\u^{\uparrow,\ell}_{0} = \mathbf{w}^{\ell}_{0}-\v^{\uparrow,\ell}_{0} $
            \EndFor
            \State $\underline{\u}^{\uparrow} =  \left (\u^{\uparrow,2}_{0} , \u^{\uparrow,2}_{1}, \u^{\uparrow,3}_{0}, ..., \u^{\uparrow,L}_{0}, \u^{\uparrow,L}_{1}  \right)^{t} $
        \EndFunction
    \end{algorithmic}
\end{algorithm}

\begin{algorithm} Upward reflections, application of $ \mathbf{\underline{L}}$ \label{alg:upwardsReflections}
    \begin{algorithmic}[1]
        \Function{ $\underline{\u}^{\uparrow}$ = Upward Reflections}{ $\underline{\v}^{\downarrow}$ }
            \For{  $\ell = 2:L-1$ }
                \State $ \begin{array}{ll} \f^{\ell} = & \delta(z_{1}-z)\v^{\downarrow,\ell}_{0} - \delta(z_{0}-z)\v^{\downarrow,\ell}_{1}  \\
                                                       & - \delta(z_{n^{\ell}+1}-z)\v^{\downarrow,\ell+1}_{0}  + \delta(z_{n^{\ell}}-z)\v^{\downarrow,\ell+1}_{1}
                                                    \end{array}$
                \State $ \mathbf{w}^{\ell}  = (\mathbf{H}^{\ell})^{-1} \f^{\ell} $
                \State $\u^{\uparrow,\ell}_{1} = \mathbf{w}^{\ell}_{1} - \v^{\downarrow,\ell}_{1}$
                \State $\u^{\uparrow,\ell}_{0} = \mathbf{w}^{\ell}_{0}  $
            \EndFor
            \State $ \f^{L} =  \delta(z_{1}-z)\v^{\uparrow,L}_{0} - \delta(z_{0}-z)\v^{\uparrow,L}_{1} $
                \State $ \mathbf{w}^{L}  = (\mathbf{H}^{L})^{-1} \f^{L} $
                \State $\u^{\uparrow,L}_{1} = \mathbf{w}^{L}_{1} - \v^{\downarrow,L}_{1}$
                \State $\u^{\uparrow,L}_{0} = \mathbf{w}^{L}_{0}  $
            \State $\underline{\u}^{\uparrow} =  \left (\u^{\uparrow,2}_{0} , \u^{\uparrow,2}_{1}, \u^{\uparrow,3}_{0}, ..., \u^{\uparrow,L}_{0}, \u^{\uparrow,L}_{1}  \right)^{t} $
        \EndFunction
    \end{algorithmic}
\end{algorithm}

\subsubsection{Physical intuition}

We deliberately present the preconditioner  in a purely algebraic fashion, as it is instructive for implementing the method.
However, there is a physical interpretation of the steps in the preconditioner, which we describe below.

As alluded to previously, the application of the preconditioner, and in particular the block back-substitution in Algs.~\ref{alg:downwards} and \ref{alg:upwards}, can be seen as a sequence of depth extrapolation steps.
Indeed, lines $4$ and $5$ in Alg.~\ref{alg:downwardsSweep} are the discrete counterpart of the incomplete Green's integral defined by
\begin{equation}
    u^{\downarrow}(\x) = \int_{\Gamma_{\ell-1, \ell}} \left (G^{\ell}(\x,\y) \partial_z  u^{\downarrow}(\y) -  \partial_z  G^{\ell}(\x,\y)  u^{\downarrow}(\y) \right) dS_y,
\end{equation}
which is equivalent to the Rayleigh integral used to extrapolate a wavefield measured in surface towards the interior of the Earth by \cite{Berkhout:Seismic_Migration_imaging_of_acoustic_energy_by_wavefield_extrapolation}.
Likewise, Lines 4 and 5 of Alg.~\ref{alg:upwardsSweep} are the discrete counterpart to an up-going discrete Green's integral.

The quality of the extrapolation depends directly on the quality of the approximation of the local Green's function $G^{\ell}$ with respect to the global Green's function.
In the reductive case, if the local Green's function is precisely the global Green's function, the method will converge in two iterations, see \cite{Gander:Optimized_Schwarz_Methods}.
However, this is equivalent to solving the global problem,  which is prohibitively expensive.
Instead, we compute a local approximation of the Green's function such that the Green's representation formula is valid within the layer only, not globally.
As expected with domain decomposition methods, incorrect local approximations introduce numerical artifacts, which are typically due to truncating the domain in a manner that is inconsistent with the underlying physics.
In the method of polarized traces, these issues are mitigated with judicious use of high-order absorbing boundary conditions in the form of PML's.
As a physical consequence, the local Green's function can only see local features within a particular layer.
Far-field interactions, reflections induced by material changes in the other layers, will not be observed by the local Green's function and must be handled iteratively, by sequentially sweeping through the domains.

An important consequence of the Green's integral representation is that it completely eliminates the difficulties that most domain decomposition methods have with seamlessly connecting sub-domains together.
Rather than assigning data dependent boundary conditions, the coupling is performed using potentials defined on the physical interfaces, and the absorbing boundary conditions in an extended domain effectively dampen spurious reflections.
The transmission conditions given by the discrete Green's representation formula are algebraically exact, thus there is no need for tuning parameters.

\section{Parallelization Strategies}

The computational effort needed to solve industrial scale 3D problems requires aggressive parallelization and optimization of the algorithm.
To obtain a scalable implementation, the algorithm and code must be designed to balance the utilization and occupancy of three key resources: CPU, memory, and communication network.
In this section, we describe our parallel implementation of the method of polarized traces, with a focus on maximizing the utilization of these resources.

In the previous section, we formally introduced a matrix-free approach for preconditioning the SIE system on layer interfaces. However, this approach still relies on local solves that are implemented using a direct solver. It is possible to use specially designed iterative local solvers by nesting the method of polarizes traces within each layer \citep{ZepedaDemanet:Nested_domain_decomposition_with_polarized_traces_for_the_2D_Helmholtz_equation} or a recursive version of the sweeping factorization \citep{Liu_Ying:Recursive_sweeping_preconditioner_for_the_3d_helmholtz_equation}. However, such approach would require a complicated code with a very carefully implemented communication pattern.
For simplicity, and to broaden the portability of the framework, we use a hybrid approach, where the local solves use off-the-shelf numerical linear algebra libraries and the polarization is matrix-free.
We will address the parallelism on two fronts: parallelism by layer and parallelism within the layers.

\subsection{Pipelining}

\begin{figure}[t]
    \begin{center}
        \includegraphics[clip, width=8.5cm]{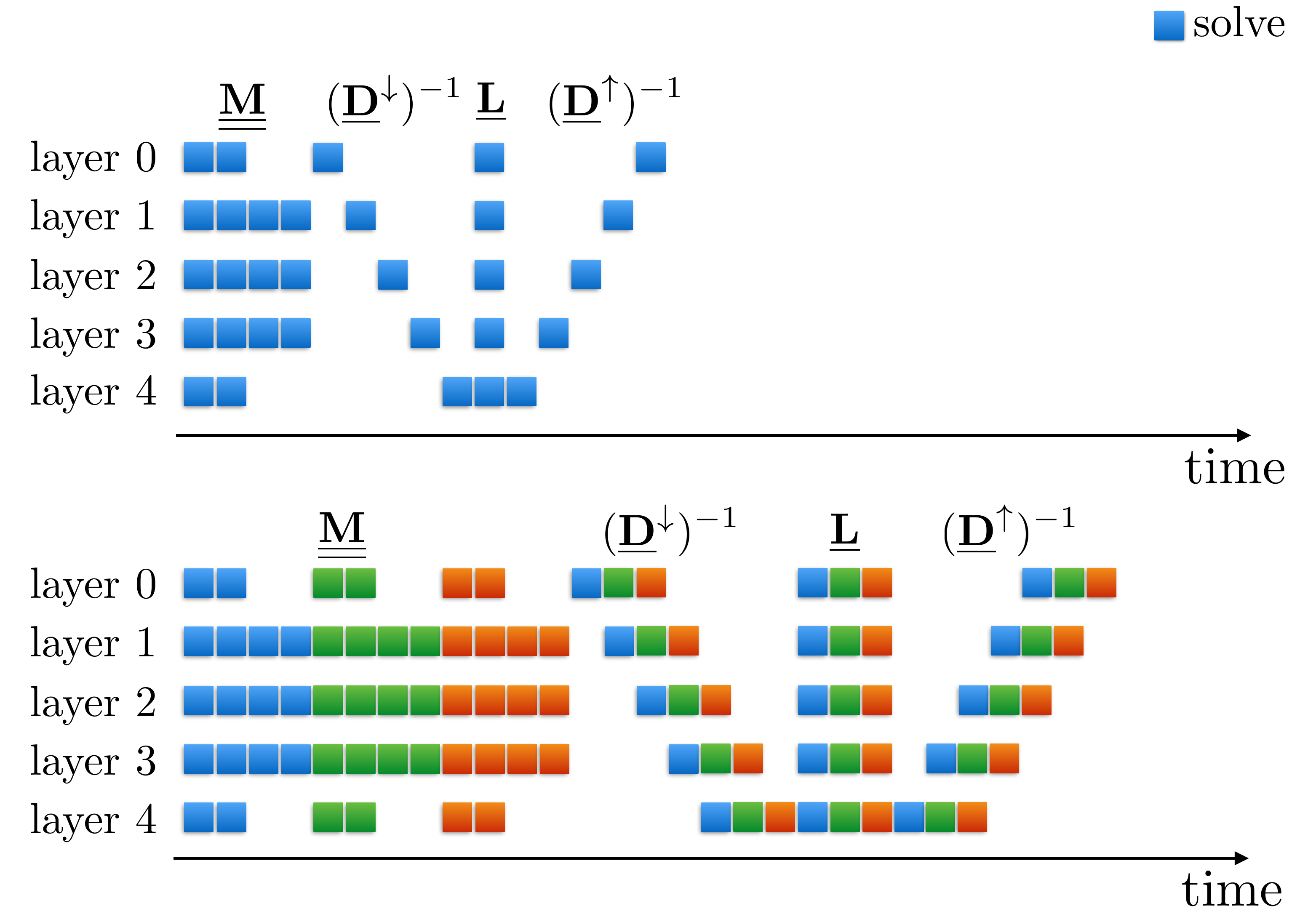}
    \end{center}
    \caption{Sketch of the load of each node in the GMRES iteration.} \label{fig:nodes_load}
\end{figure}

First, we address parallelism due to the layer decomposition.
Primarily, the parallelism across layers is due to the SIE and the preconditioner used to help solve it.
There are five trivially parallel (by layer) applications of the local solver: four due to $\underline{\underline{\mathbf{M}}}$ and one due to the appearance of $\underline{\mathbf{L}}$ in the preconditioner.
However, in the preconditioner application there are sequential bottlenecks due to the applications of  $(\mathbf{\underline{D}}^{\uparrow})^{-1}$ and $(\mathbf{\underline{D}}^{\downarrow})^{-1}$ via block back-substitution.
Despite the trivial parallel nature of the other local solver applications, applying the preconditioner using Algs.~\ref{alg:downwardsSweep} and \ref{alg:upwardsSweep} permits work to be done on only one layer at a time,
thus forcing the majority of the computer to remain idle.
This is illustrated in the top half of Fig.~\ref{fig:nodes_load}, where each blue box represents a local solve and algorithm execution moves from left to right.
Supposing that each local solve costs $\gamma(n)$ time, then following Fig.~\ref{fig:nodes_load}, each GMRES iteration can be performed in $5 \gamma(n) + 2L \gamma(n)$, ignoring communication costs.

To alleviate the sequential bottleneck, we leverage the fact that seismic problems have thousands of right-hand sides and introduce pipelining.
Pipelining allows us to process multiple right-hand sides simultaneously, each at different levels of progress through the sweeps, which helps balance the computational load on the layers, reducing the idle time and increasing the computational efficiency.
The pipelining principle is demonstrated in the bottom half of Fig~\ref{fig:nodes_load}, where the boxes represent a local solve and the blue, green, and orange colors indicate distinct right-hand sides.
Pipelining allows the layers to perform work for different right-hand sides simultaneously.
Indeed, as long as there are at least $2L$ right-hand sides, the pipeline can be completely full and all available compute resources will be occupied.
For the pipelined algorithm, again disregarding communication costs, the runtime of a GMRES iteration is $5R\gamma(n) + 2(L + R)\gamma(n)$.
Recalling that $L \sim n$, $R\sim n$, and $\gamma(n) = \cO(\alpha^2_{\text{pml}} n^2 \log n) $, the cost ratio for solving $R$ right-hand sides compared to one right-hand side is constant,
\begin{equation}
   \frac{5R\gamma(n) + 2(L + R)\gamma(n)}{5 \gamma(n) + 2L \gamma(n)} = \cO(1).
\end{equation}

One of the advantages of the method of polarized traces is that the memory requirement to store the intermediate representation of solution is lower than other methods, because it requires solutions for the degrees of freedom involved in the SIE only.
Thus, for each right-hand side, only $N/q$ data need to be stored, where $q$ is the thickness of the interface.
This reduction in storage, when combined with the reduced storage due to the relatively small number of GMRES iterations required for convergence, yields a smaller memory footprint for the outer GMRES iteration than methods requiring to update the full volume.
It is possible to further reduce the memory footprint by using Bi-CGSTAB instead of GMRES, keeping the computational cost almost constant \citep{ZepedaDemanet:Nested_domain_decomposition_with_polarized_traces_for_the_2D_Helmholtz_equation}.

\subsection{Local Solves: Parallel Multi-frontal Methods}

The method polarized traces is a highly modular framework for solving the Helmholtz equation: in practice, one can use any existing algorithm or package for solving linear equations to perform the local solves, including the method of polarized traces itself~\citep{ZepedaDemanet:Nested_domain_decomposition_with_polarized_traces_for_the_2D_Helmholtz_equation}.
To obtain good parallel performance, we use a high-performance distributed linear algebra library to solve the local problems within each layer.
Due to the sparsity pattern of the linear system at each layer, a typical recipe for the local solves is to re-order the degrees of freedom to increase stability and reduce numerical fill-in, perform a multi-frontal factorization, and solve the resulting factorized system with forward- and backward-substitution, often called triangular solves.
There exists a myriad of techniques to parallelize the multi-frontal factorization and the triangular solves (see \cite{davis:A_survey_of_direct_methods_for_sparse_linear_systems} for a recent and extensive review of different techniques for solving sparse systems).
A popular approach is to use supernodal elimination trees (see \cite{Ashcraft:Progress_in_Sparse_Matrix_Methods_for_Large_Linear_Systems_On_Vector_Supercomputers}) defined through nested dissection, which results in highly scalable factorizations \citep{Gupta:Highly_Scalable_Parallel_Algorithms_for_Sparse_Matrix_Factorization}, albeit with less efficient triangular solves \cite{Gupta_Kumar:A_high_performance_two_dimensional_scalable_parallel_algorithm_for_solving_sparse_triangular_systems}.
To avoid the poor scalability of dense triangular solves due to this approach, \citet{Raghavan:Efficient_Parallel_Sparse_Triangular_Solution_Using_Selective_Inversion} introduced a scheme called {\it selective inversion}, which is applied by \citet{Poulson_Engquist:a_parallel_sweeping_preconditioner_for_heteregeneous_3d_helmholtz} specifically for the Helmholtz problem. Although very efficient, using the techniques mentioned before would require lengthy and complex code, we use instead off-the-shelf libraries, which can be effortlessly changed if necessary.

For the results in this paper, we use STRUMPACK \citep{Rouet_Li_Ghysels:A_distributed-memory_package_for_dense_Hierarchically_Semi-Separable_matrix_computations_using_randomization} to perform the local solves.
STRUMPACK is a state-of-the-art distributed sparse linear solver library that relies on supernodal factorizations,\footnote{ It uses a $ULV$ factorization when the compression is turned on \citep{Chandrasekaran:A_Fast_ULV_Decomposition_Solver_for_Hierarchically_Semiseparable_Representations,Xia:Randomized_Sparse_Direct_Solvers}}
a 2D block-cyclic distribution of the matrix, and a static mapping technique to assign task to MPI processes based on {\it proportional mapping} \citep{Pothen_Sun:A_Mapping_Algorithm_for_Parallel_Sparse_Cholesky_Factorization} to achieve good parallel performance.
The implementation is competitive with other distributed linear algebra solvers with liberal licenses, such as SuperLU-DIST (cf. \cite{li_demmel03:SuperLU_DIST}) and MUMPS (cf. \cite{Amestoy_Duff:MUMPS}), while providing the user more freedom to arrange the distribution of the matrix and right-hand sides, within a distributed memory enviroment, in a manner that is optimal for specific applications.

A strong advantage of STRUMPACK is that the factorization and solve processes can be accelerated using compressed linear algebra, in particular HSS compression with nested bases, using randomized sampling techniques.
In general, using compressed formulations allows a reduced memory footprint and in some cases faster algorithms.
For the high-frequency regime, it is known that solvers based on compression techniques do not provide a lower asymptotic complexity, due to the fact that the ranks of the off-diagonal blocks are frequency dependent \citep{Engquist_Zhao:approximate_separability_of_green_function_for_high_frequency_Helmholtz_equations}.
However, these solvers still tend to provide smaller memory footprints for the cases we consider in this paper, albeit, with much bigger runtimes.
Therefore, we deliberately do not considered the performance of the adaptive compression in this paper and leave such treatment for future work.

The STRUMPACK solver addresses parallelism within layers in two ways: classical distribution of tasks with MPI and synchronous processing within each task with OpenMP, which we exploit for the results presented in our numerical results.
STRUMPACK's hybrid parallelism model allows us to maximize the utilization of computational resources.
We have designed the distribution of MPI tasks to exploit highly asynchronous communication patterns, thus reducing the communication time subtantially.

Finally, in most of the experiments shown in the sequel, we process at most one right-hand side per layer. It is possible to solve more than one right-hand side per layer, and take advantage of BLAS3 routines. Moreover, a quick inspection to the algorithms  within the preconditioner, shows that the right-hand sides are sparse. Indeed, the sources are supported on the interfaces, and the solutions are only needed at the boundaries. In principle, it is possible to take advantage of the sparsity of the solution and the right-hand-side to reduce the constants by removing some branches from the elimination tree. These techniques would certainly reduce the constants but they would have little impact on the asymptotic scaling; hence, they were not explored in this study.

\subsection{Communication}

\begin{figure}[t]
    \begin{center}
        \includegraphics[trim = 35mm 245mm 110mm 10mm, clip, width=8.5cm]{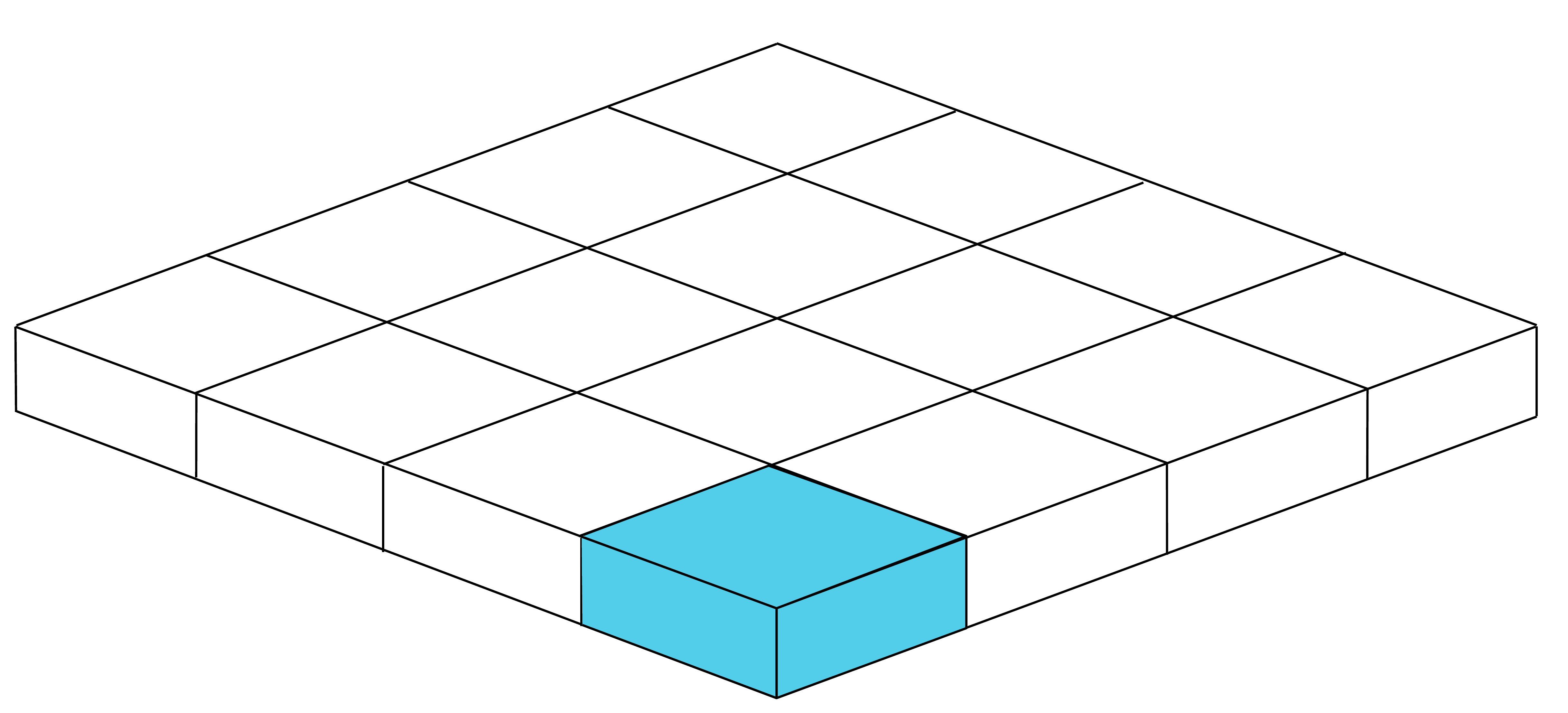}
    \end{center}
    \caption{Sketch of the decomposition of the degrees of freedom of the slabs in cubes.} \label{fig:DDM_cubes}
\end{figure}

As is common in massively parallel applications, there is a bottleneck due to the communication between parallel tasks, which is strongly dependent on the distribution of the tasks on the cluster.
For this implementation, we assume a very simple topology for the distribution of unknowns.
We distribute the parallel tasks following the layer structure.
For each slab of unknowns, we assign $\cO(n^2)$ tasks, and using MPI directives enforce that the tasks are contiguous within physical computing nodes.
Each slab is divided into $\cO(n^2)$ cubes, as illustrated in Fig.~\ref{fig:DDM_cubes}, and the parallel tasks associated with that slab are divided evenly and contiguously amongst the cubes.
Each cube contains a contiguous block of the solution, as shown in Fig.~\ref{fig:DDM_cubes}, and the associated entries of the local matrix.

\begin{figure}[t]
    \begin{center}
        \includegraphics[trim = 50mm 35mm 27mm 18mm, clip, width=8.5cm]{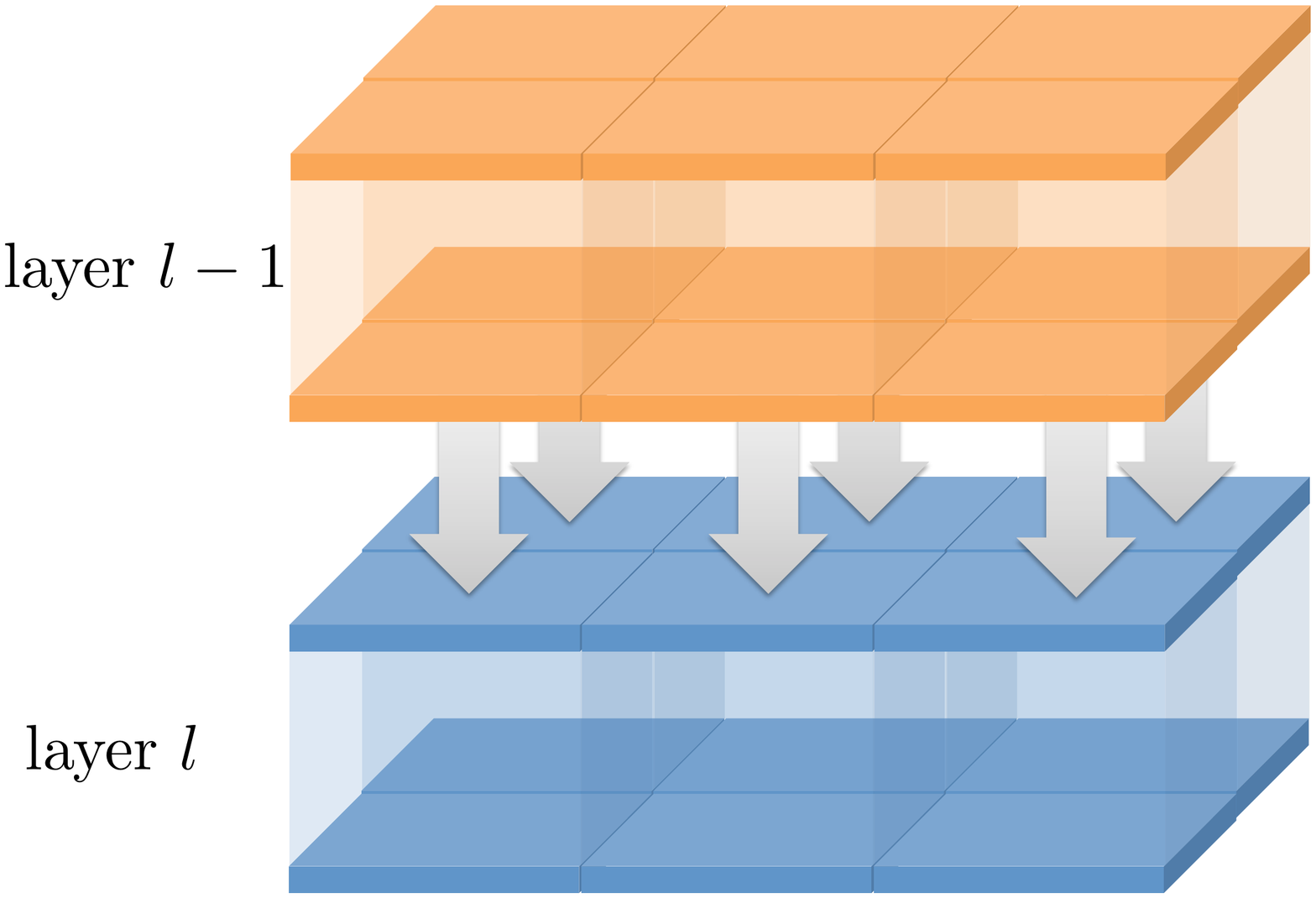}
    \end{center}
    \caption{Sketch of the asynchronous communication between slabs.} \label{fig:communication}
\end{figure}

Within the slab, the topology is designed so that each cube only communicates with its neighboring cubes in the same slab, generally inside the local solver.
Across slabs, cubes only communicate with the cube in the same position on the adjacent slabs immediately above and below it, as illustrated in Fig.~\ref{fig:communication}.
Under this particular topology, we can distinguish two main communications bottlenecks:
\begin{itemize}
\item the communication between parallel tasks within the distributed linear algebra solver, and
\item the communication of the boundary data between slabs, during the application of the preconditioner.
\end{itemize}

As a consequence of using third-party distributed linear algebra solvers, we have little control over the communication pattern, particularly because STRUMPACK uses MC64 \citep{Duff_Koster:On_Algorithms_For_Permuting_Large_Entries_to_the_Diagonal_of_a_Sparse_Matrix}, for enhancing stability and ParMetis \citep{Karypis:A_Parallel_Algorithm_for_Multilevel_Graph_Partitioning_and_Sparse_Matrix_Ordering} to optimally reorder the matrix in order to reduce fill-in during the factorization.
To have the desired distribution of the degrees of freedom among the cubes, we reorder the matrix with a $Z$ ordering, so that smallest division corresponds exactly to the degrees of freedom within a cube.
Then, the matrix is then assembled and passed to the linear solver in a distributed fashion.

Communication between slabs is a product of the application of the preconditioner in Eq.~\ref{eq:preconditioner_GS}, in which back- and forward-substitution are used to apply $(\mathbf{\underline{D}}^{\uparrow})^{-1}$ and $(\mathbf{\underline{D}}^{\downarrow})^{-1}$.
Algs.~\ref{alg:downwards} and \ref{alg:upwards}) require a local solve in each slab, followed by communication of the trace information to the next slab, in which another local solve is performed. This operation is repeated until all the slabs are visited within the sweep.

By dividing the slabs into cubes, the communication of the trace information between slabs is very efficient.
Given that each cube communicates with the cube directly above and below, it is possible to perform asynchronous point-to-point communication between the cubes of two adjacent slabs as shown in Fig.~\ref{fig:communication}.
This allows the communication to be performed in nearly constant time, up to saturation of the network.
Moreover, the trace information is already distributed for distributed assembly of the right-hand side within the subsequent slab.
As stated before, it is possible to use topologies better suited for the multi-frontal solver, such as the one due to \citet{Poulson_Engquist:a_parallel_sweeping_preconditioner_for_heteregeneous_3d_helmholtz}.
However, such an implementation requires a very precise understanding of the re-ordering mechanism within the solver, which we do not have for black-box solvers.

\section{Complexity of Polarized Traces}

The run-time complexity of the polarized traces preconditioner is driven by the costs of computation and communication.
Achieving optimal performance requires delicately balancing the parallel distribution of the problem depending on the characteristics of the target HPC system.
In this section we develop models for computation and communication costs, which guide problem parameter selection in HPC environments.

\subsection{Computation}

As before, each layer has $\cO((n_z + \alpha_{\text{pml}})\times n^2)$ grid points, i.e., they are $n_z$ grid points thick with $\alpha_{\text{pml}}$ additional points due to the PML.
We have that $n_z = \cO(1)$ because $L \sim n$, which implies that we are solving a quasi-2D problem.
The additional cost is due only to the points used to implement the absorbing boundary conditions. When applying 2D nested dissection to the quasi-2D problem we have $\cO(\alpha_{\text{pml}} n)$ degrees of freedom in the biggest front, thus leading to a complexity of $\cO(\alpha_{\text{pml}}^3n^3)$ for the factorization of the systems local to each layer and $\cO(\alpha_{\text{pml}}^2n^2 \log{n})$ for the application of the triangular solve \citep{Duff_Reid:The_Multifrontal_Solution_of_Indefinite_Sparse_Symmetric_Linear}.
Sequentially, the complexity of Alg.~\ref{alg:offline_computation} is $\cO(\alpha_{\text{pml}}^3 N^{4/3})$, but given that the loop in line $2$ of Alg.~\ref{alg:offline_computation} is embarrassingly parallelizable, Alg.~\ref{alg:offline_computation} can be performed in $\cO(\alpha_{\text{pml}}^3 N)$ time\footnote{As it will be shown in the numerical experiments, this scaling can be further reduced due to the parallelism at the level of the multi-frontal solver.}.
Due to the sequential nature of Algs.~\ref{alg:downwardsSweep} and~\ref{alg:upwardsSweep}, applying the preconditioner requires $2L$ local solves per iteration, applied sequentially.
Consequently, the total complexity for the application of the preconditioner is $\cO(\alpha_{\text{pml}}^2 N \log{N})$.
For $\alpha_{\text{pml}} \sim \log{n} $, at most $\cO(\log{n})$ iterations are empirically needed to converge, thus the complexity of the solver is linear (up to poly-logarithmic factors), provided that $L \sim n$ and that the number of iterations for convergence grows slowly.

It is possible to relax the restriction that $L \sim n$, instead allowing $L \sim n^{b}$ where $b < 1 $.
However, in this regime, maintaining the overall linear complexity requires that we exchange the multi-frontal solver for an iterative solver \citep{Liu_Ying:Recursive_sweeping_preconditioner_for_the_3d_helmholtz_equation,ZepedaDemanet:Nested_domain_decomposition_with_polarized_traces_for_the_2D_Helmholtz_equation}.
The main disadvantage of this approach is that it reduces the possible parallelism due to using multi-frontal solvers, which makes an efficient implementation of the pipelining difficult and makes the communication patterns more complicated.

\subsection{Pipelining}

We introduced pipelining of $R$ right-hand sides to alleviate the sequential nature of applying the preconditioner.
To understand the run-time impact of pipelining multiple right-hand sides, we consider its impact on the complexity of applying $\underline{\underline{\mathbf{M}}}$ and applying the preconditioner.
Given that the run-time cost of each local solve is $\cO(\alpha_{\text{pml}}^2n^2 \log n)$ and recalling that applying $\underline{\underline{\mathbf{M}}}$ is embarrassingly parallel, the cost of applying $\underline{\underline{\mathbf{M}}}$ to $R$ right-hand sides is $\cO(\alpha_{\text{pml}}^2 R n^2 \log{n})$.
As long as $R \sim L \sim n $, applying the preconditioner costs $\cO(L \alpha_{\text{pml}}^2 n^2 \log {n})$.
However, when $R \gtrsim L$, the additional right-hand-sides are treated sequentially, resulting in a cost of $\cO(\alpha_{\text{pml}}^2 R n^2 \log {n})$.
Using the fact that $L \sim n$ and that $N = n^3$, we obtain the advertised runtime of $\cO( \alpha_{\text{pml}}^2  \max(1,R/L) N \log {N} ))$.

\subsection{Communication}

\setcounter{footnote}{0}

We treat the distributed linear algebra solver as a black box, and therefore do not analyze the costs of communication within the local solve.
Thus, we only consider the cost of communication due to the global solve, that is the costs of communicating between subdomains across layers.
We assume that that each subdomain has fast access to its corresponding patch of both the $R$ wavespeeds and the $R$ sources.
Moreover, we assume that each subdomain assembles and stores its portion of the global solution\footnote{As a consequence, computation of imaging conditions can be performed without incurring in extra communication cost.}.

The offline stage of the algorithm, assembly of the local matrices and the local factorizations, is embarrassingly parallel under the assumptions described above.
The online part, has three stages:
\begin{itemize}
\item  the preparation of the right-hand side  (Lines 2-8 in Alg.~\ref{alg:SIE_solver});
\item  the solve of the SIE (Lines 9 in Alg.~\ref{alg:SIE_solver}); and
\item  the assembly of the global solution (Lines 10-14 in Alg.~\ref{alg:SIE_solver}).
\end{itemize}
Of these, the first and third stages require no communication under the above assumptions.

Solving the SIE has two main phases: the application of $\underline{\underline{\mathbf{M}}}$ and the application of the preconditioner.
Applying $\underline{\underline{\mathbf{M}}}$, as shown in Alg.~\ref{alg:applyM}, is an embarassingly parallel operation with zero communication.
On the other hand, applying the preconditioner, which is fully sequential, requires communication of $\cO(n^2)$ unknowns from one layer to the next.
In our implementation, these unknowns are distributed evenly between $\cO(n^2)$ MPI tasks.
Using a point-to-point communication strategy, each of the MPI tasks assigned to a layer only communicates with one corresponding MPI task in the adjacent (above and below) layers, as illustrated in Fig.~\ref{fig:communication}.
Thus, by exploiting asynchronous communication, the communication between layers can be performed in $\cO(1)$ time up to saturation of the bandwidth, which is asymptotically negligible with respect to solving the local linear systems.
This communication must be performed $\cO(L)$ times during each sweeping operation in Algs.~\ref{alg:downwardsSweep} and~\ref{alg:upwardsSweep}.
Consequently, total communication cost of applying the preconditioner is $\cO(n)$, up to saturation of the bandwidth.

\section{Numerical Experiments}

In this section, we present the results of several numerical experiments used to verify the complexity described above.
In particular, we demonstrate the performance of the 3D preconditioner in various heterogeneous media for a single source and then illustrate the impact of pipelining on parallel performance.
Our polarized traces implementation is written in C and compiled with the 2015 Intel compiler suite.
The current implementation is uses IEEE double precision floating point.
To perform the local solves, we use STRUMPACK v1.1.0 with Intel MKL support for fast linear algebra operations.
The preconditioner is parallelized with MPI and STRUMPACK is parallelized with both MPI and OpenMP.
The experiments were performed on Total's ``Laure'' SGI ICE-X cluster, where each computing node contains two 8-core Intel Sandy Bridge processors, 64GB of RAM, and are connected with an Infiniband interconnect.

\subsection{Homogeneous media}

First, we demonstrate the effectiveness of the preconditioner by solving the Helmholtz problem in homogeneous media. With no reflectors in the medium, the convergence of the algorithm is only dependent on the frequency and the quality of the absorbing boundary condition at the layer interfaces.
In this experiment, as well as the subsequent experiments, we hold $\alpha_{\text{pml}}$ constant, with sufficient points to minimize artificial reverberations while simultaneously keeping the number of iterations low.
For higher frequencies, to preserve the low iteration count we would need to scale $\alpha_{\text{pml}}$ as $\cO(\log{n})$.

In this experiment, we test 4 problem sizes, $n=50,$ $100,$ $200,$ and $400$, which corresponds to frequencies of $8$, $16$, $32$, and $64$ Hz.
Source frequency is scaled with problem size to stay in the high-frequency regime and sources are assumed to be point sources.
The number of layers is also scaled with the problem, $L=5,$ $10,$ $20,$ and $40$.
Due to memory limitations on the computing node, in some cases the nodes were saturated before all cores could be assigned to an MPI task.
For these cases, we allow the remaining cores to be used for vectorized processing with OpenMP.
The outer GMRES iteration is run until the residual is reduced to $10^{-7}$, which is excessive in a production, single-precision environment; however, it shows the favorable behavior of the solver under more challenging conditions. Lower tolerances can produce misleading results when frequencies are not high enough, thus only revealing a pre-asymptotic behavior.
For each configuration,  we report the wall-clock times for initialization, matrix assembly, matrix factorization, and total online time for $R=1$ and $R=L$ with pipelining.
Additionally, we track the number of GMRES iterations required to achieve the desired convergence.

\begin{figure}[t]
    \begin{center}
        \includegraphics[trim = 0mm 10mm 0mm 10mm, clip, width=7.5cm]{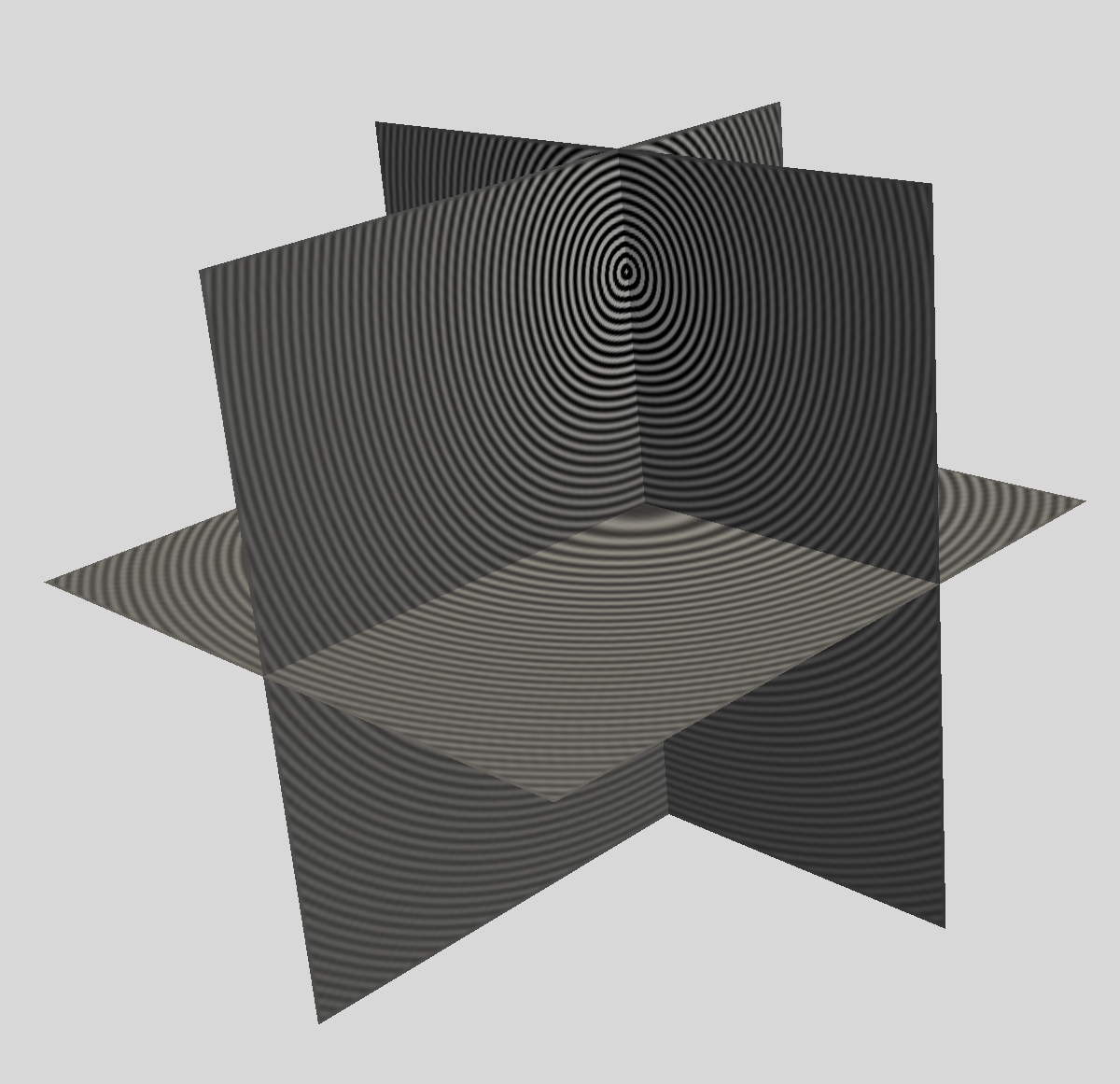}
    \end{center}
    \caption{Solution of the Helmholtz equation, at 64 Hz, in constant wavespeed.} \label{fig:constant_model}
\end{figure}

\begin{figure}[t]
    \begin{center}
        \includegraphics[clip,width=8.5cm]{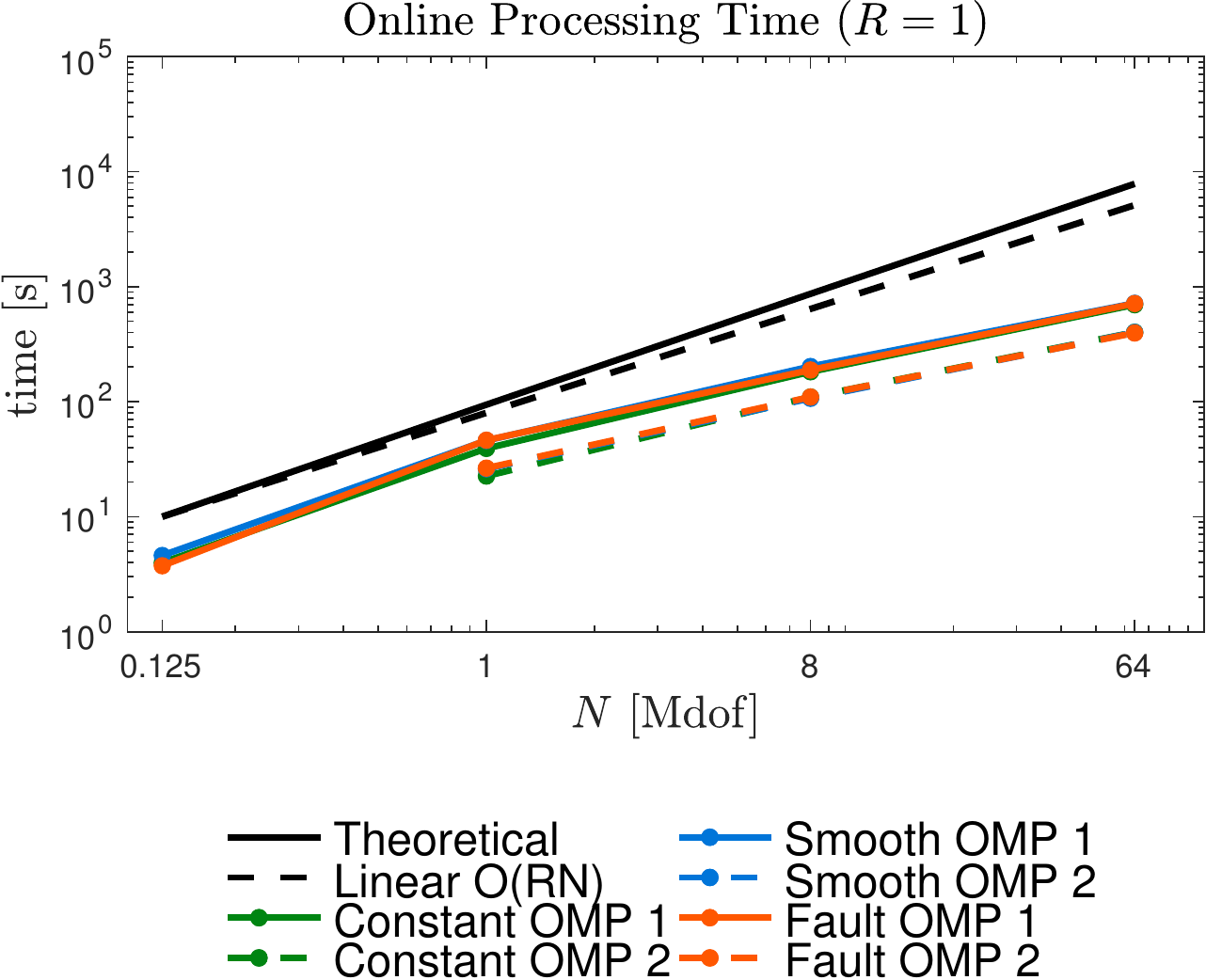}
    \end{center}
    \caption{Observed run-time as function of $N$ for homogeneous media (green),  smooth heterogeneous media (blue) and ``fault'' model  (orange), with pure MPI (solid) and hybrid MPI-OpenMP (dashed) for $R=1$ right-hand sides.  For comparison, theoretical scaling of polarized traces algorithm is given (solid black) as well as linear scaling (dashed black).} \label{fig:not_seam_1}
\end{figure}

\begin{figure}[t]
    \begin{center}
        \includegraphics[clip, width=8.5cm]{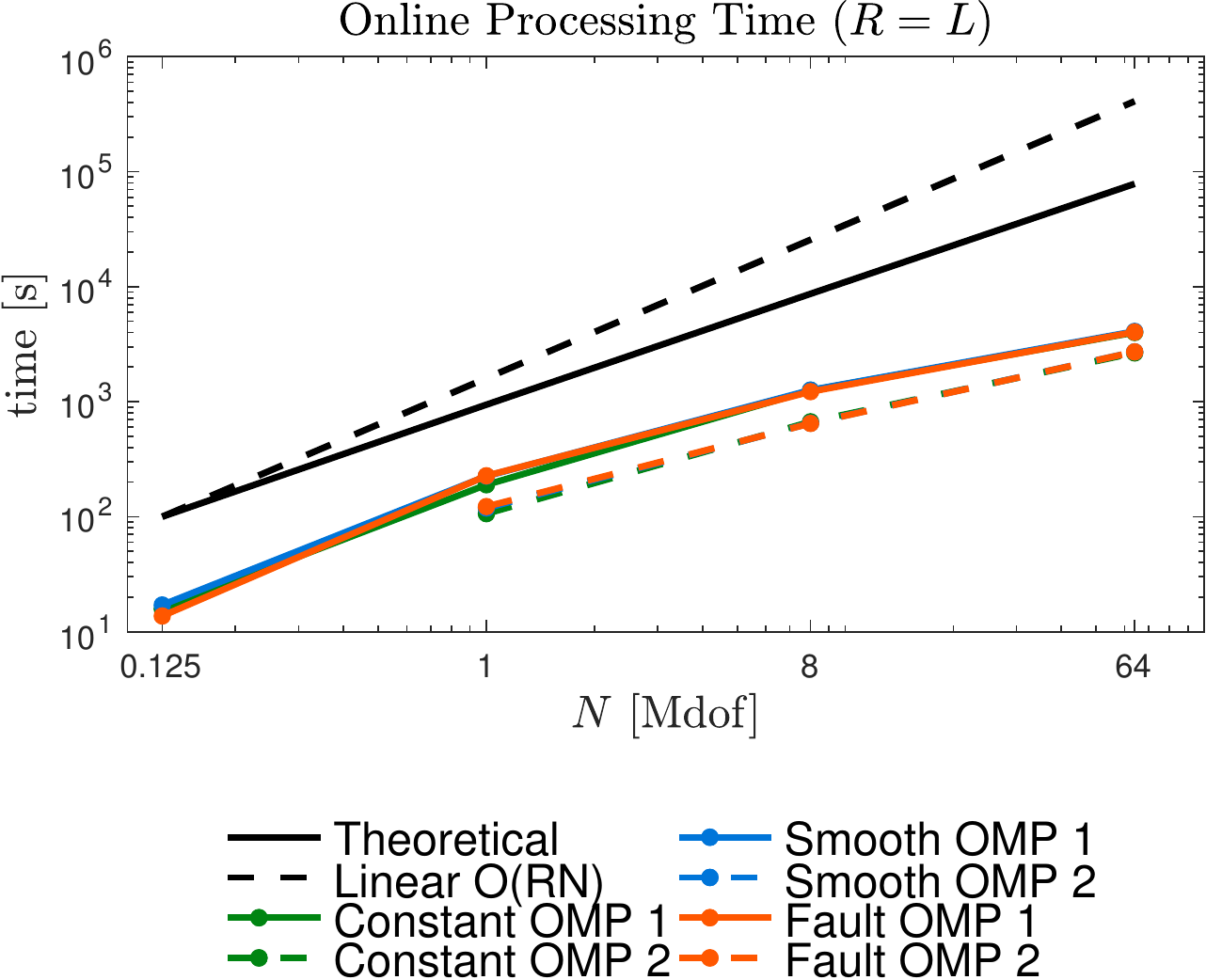}
    \end{center}
    \caption{Observed run-time as function of $N$ for homogeneous media (green),  smooth heterogeneous media (blue) and ``fault'' model  (orange), with pure MPI (solid) and hybrid MPI-OpenMP (dashed) for $R=L$ right-hand sides.  For comparison, theoretical scaling of polarized traces algorithm is given (solid black) as well as linear scaling (dashed black).} \label{fig:not_seam_R}
\end{figure}

\begin{table*}[t]
    \begin{center}
\begin{tabular}{|r|cccccccc|}
\hline
N ~      & $50^3$ &  $100^3$ &  $100^3$ &  $200^3$ & $200^3$ &  $400^3$ & $400^3$  & $400^3$ \\
\hline
L  ~    & 5 & 10 & 10 & 20 & 20 & 40 & 40 & 40 \\
\hline
MPI Tasks   & 5 & 10 & 10 & 80 & 80 & 640 & 640 & 640\\
OpenMP Threads per Task     & 1 & 1 & 2 & 1 & 2 & 1 & 2 & 3\\
Total Cores & 5 & 10 & 20 & 80 & 160 & 640 & 1280 & 1920\\
Total Nodes & 1 & 1 & 2 & 5 & 10 & 80 & 80 & 128\\
\hline
\textbf{Single rhs} &  &  &  &  &  &  &  & \\
\hline
\# GMRES Iterations & 4 & 4 & 4 & 5 & 5 & 6 & 6 & 6\\
 Initialization [s]         & 0.2 & 1.0 & 0.9 & 6.9 & 4.4 & 18.9 & 18.9 & 18.4\\
 Factorization [s]      & 4.1 & 41.1 & 21.9 & 153.2 & 78.3 & 320.5 & 200.1 & 148.6\\
 Online [s]                 & 4.0 & 39.2 & 22.6 & 182.0 & 109.7 & 696.6 & 401.4 & 315.5\\
 Average GMRES [s]  & 0.9 & 8.4 & 4.8 & 32.0 & 19.2 & 103.5 & 59.3 & 46.6\\
 \hline
 \textbf{Pipelined rhs} &  &  &  &  &  &  &  & \\
\hline
 $R$ (number of rhs)    & 5 & 10 & 10 & 20 & 20 & 40 & 40 & 40\\
 Online [s]                     & 15.8 & 189.4 & 106.2 & 1255.5 & 668.5 & 3994.2 & 2654.4 & 1878.1\\
 Average GMRES [s]  & 3.4 & 40.6 & 22.7 & 223.8 & 118.6 & 599.9 & 401.0 & 283.0\\
 Online per rhs [s]                 & 3.2 & 18.9 & 10.6 & 62.8 & 33.4 & 99.9 & 66.4 & 47.0\\
 Average GMRES per rhs [s]  & 0.7 & 4.1 & 2.3 & 11.2 & 5.9 & 15.0 & 10.0 & 7.1 \\
 \hline
\end{tabular}
\end{center}
    \caption{Runtime (in seconds) for one and several right-hand sides for the solution of the Helmholtz equation using an homogeneous model.}\label{table:test_homogeneous}
\end{table*}

The solution at 64 Hz is provided in Fig.~\ref{fig:constant_model}.
The full results of the experiment are given in in Table~\ref{table:test_homogeneous} and the observed run-times, compared to the theoretical run-times for $R=1$ and $R=L$ pipelined right-hand sides are shown in Fig.~\ref{fig:not_seam_1} and Fig.~\ref{fig:not_seam_R}, respectively.
Only when $R \ge L$ does the theoretical scalability break the linear threshold, however, in both cases, the method of polarized traces scales better than the theoretical scaling, which we attribute to optimizations and parallelism in the local solver.
Of note in Table~\ref{table:test_homogeneous}, the number of GMRES iterations grows very slowly with the frequency, even when we do not scale $\alpha_{\text{pml}}$ optimally.
In the experiments where hybrid parallelism (MPI-OpenMP) is used, the run-times are reduced almost linearly for medium-sized problems, however the improvements fade as the size of the problem increases.
This behavior is due to the linear solver, where for large problems the memory access time in the  triangular solves becomes dominant, reducing the parallelism.

\subsection{Smooth heterogeneous media}

\begin{figure}[htb]
    \begin{center}
        \includegraphics[trim = 40mm 50mm 10mm 60mm, clip, width=8.5cm]{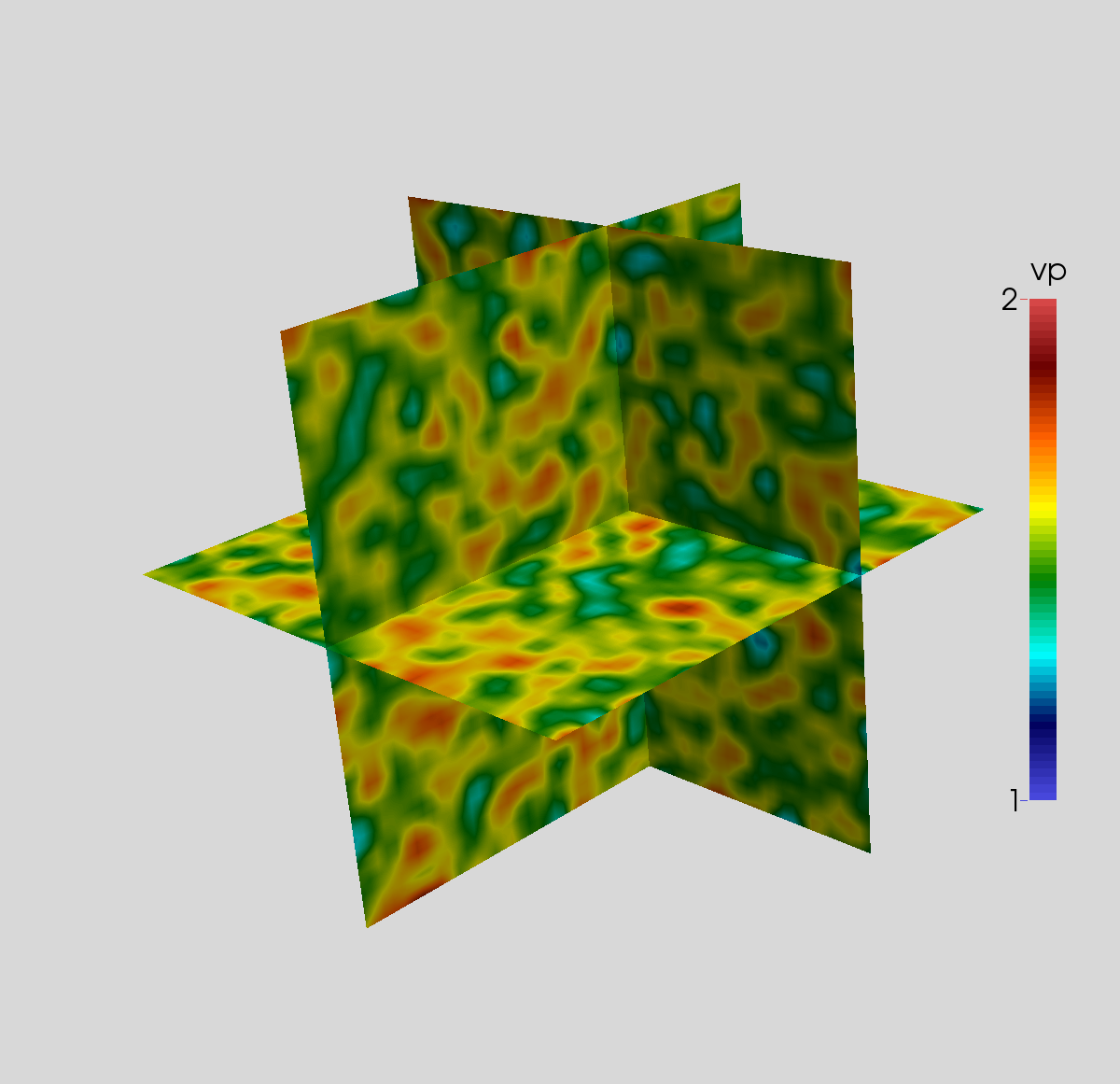}
    \end{center}
    \caption{Randomly generated smooth heterogeneous medium.} \label{fig:smooth_model}
\end{figure}

\begin{figure}[htb]
    \begin{center}
        \includegraphics[trim = 0mm 10mm 0mm 10mm, clip, width=8.5cm]{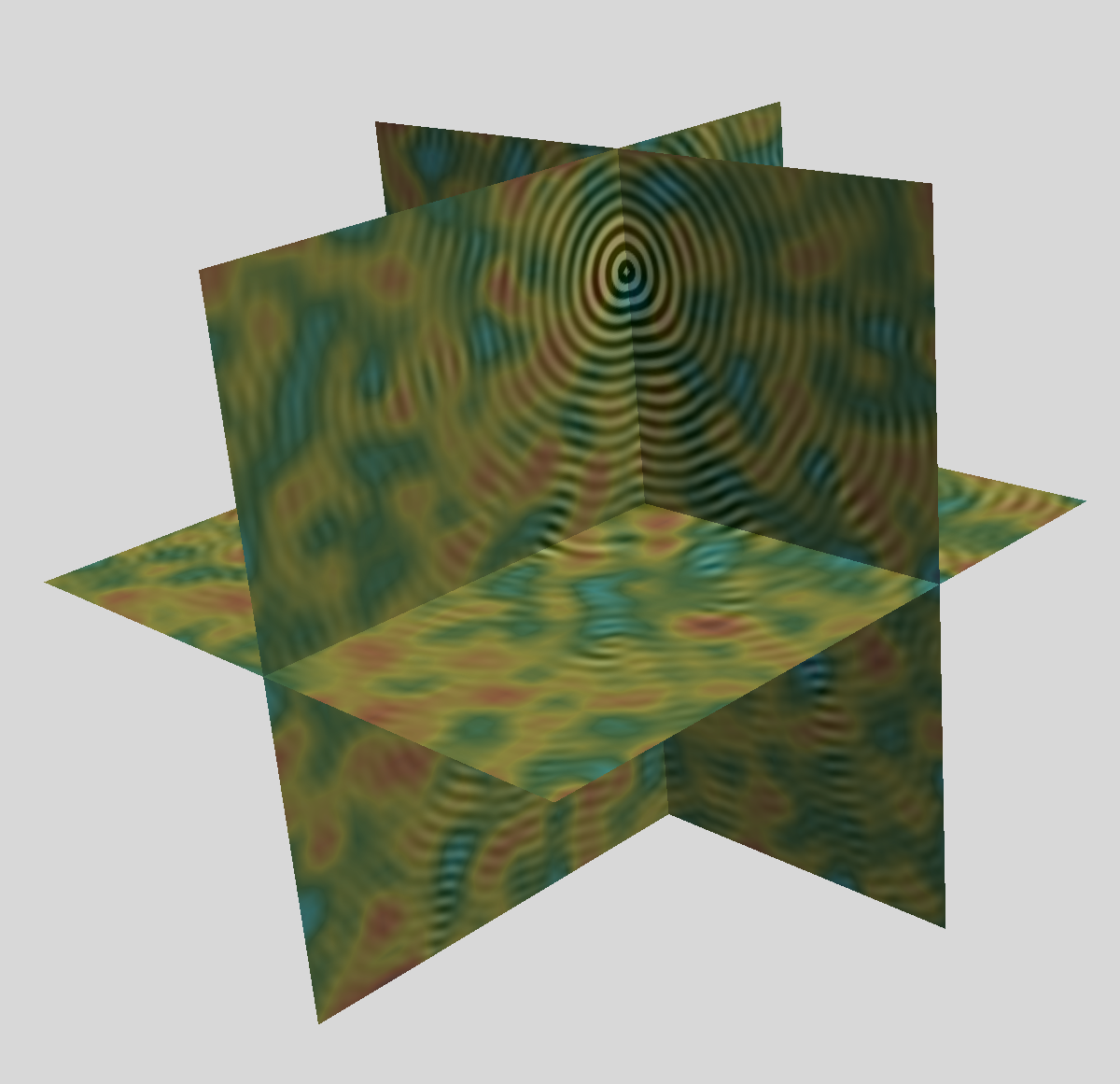}
    \end{center}
    \caption{Real part of the solution to the Helmholtz equation using the random smooth medium.} \label{fig:smooth_model_solution}
\end{figure}

Using the same configurations as above we solve the Helmholtz problem in the smoothed random media shown in Fig.~\ref{fig:smooth_model} in order to demonstrate the effectiveness of the solver in heterogeneous media.
This test demonstrates that the method is particularly robust to media where the rays can bend and develop caustics.
The solution for the configuration equivalent to that of Fig.~\ref{fig:constant_model} is given in Fig.~\ref{fig:smooth_model_solution}.
In Fig.~\ref{fig:smooth_model_solution}, it is clear that the features of the model are comparable to the wavelength used, thus the solution presents interference, caustics, non-spherical wavefronts.
Table \ref{table:test_smooth} contains the complete experimental results, where we observe that the variation in the media has little real effect on the run-time or convergence properties.
In this case, even if the rays bend, the preconditioner does a remarkable job at tracking the rays in the correct direction and propagating them accordingly.
The timings are given in in Figs.~\ref{fig:not_seam_1} and Fig~\ref{fig:not_seam_R}.

\begin{table*}
\begin{center}
\begin{tabular}{|r|cccccccc|}
\hline
N ~      & $50^3$ &  $100^3$ &  $100^3$ &  $200^3$ & $200^3$ &  $400^3$ & $400^3$  & $400^3$ \\
\hline
L  ~    & 5 & 10 & 10 & 20 & 20 & 40 & 40 & 40 \\
\hline
MPI Tasks               & 5 & 10 & 10 & 80 & 80 & 640 & 640 & 640\\
OpenMP Threads per Task & 1 & 1 & 2 & 1 & 2 & 1 & 2 & 3\\
Total Cores             & 5 & 10 & 20 & 80 & 160 & 640 & 1280 & 1920\\
Total Nodes             & 1 & 1 & 2 & 5 & 10 & 80 & 80 & 128\\
\hline
\textbf{Single rhs} &  &  &  &  &  &  &  & \\
\hline
\# GMRES Iterations     & 5     & 5     & 5     & 5     & 5     & 6     & 6     & 6\\
Initialization [s]          & 0.2   & 1.1   & 1.0   & 7.3   & 4.6   & 21.3  & 21.2  & 20.8\\
Factorization  [s]      & 3.8   & 41.1  & 21.8  & 156.0 & 79.4  & 323.7 & 204.5 & 151.5\\
Online [s]                      & 4.6   & 45.9  & 26.1  & 202.2 & 106.9 & 717.0 & 400.1 & 314.5\\
Average GMRES [s]   & 0.8   & 8.1   & 4.6   & 35.5  & 18.7  & 106.4 & 59.2  & 46.5\\
\hline
\textbf{Pipelined rhs} &  &  &  &  &  &  &  & \\
\hline
$R$ (number of rhs) & 5 & 10 & 10 & 20 & 20 & 40 & 40 & 40\\
Online [s]                      & 17.1  & 225.1 & 118.8 & 1260.9  & 650.2   & 4085.0    & 2714.8 & 1872.1\\
Average GMRES [s]   & 3.0   & 39.8  & 20.9  & 223.6   & 115.6   & 613.3     & 409.2 & 281.9\\
Online  per rhs [s]         & 3.4   & 22.5  & 11.9  & 63.0    & 32.5    & 102.1     & 67.9  & 46.8\\
Average GMRES  per rhs [s]  & 0.6   & 4.0   & 2.1   & 11.2    & 5.8     & 15.3  & 10.2  & 7.0\\
\hline
\end{tabular}
\end{center}
    \caption{Runtime (in seconds) for one and several right-hand sides for the solution of the Helmholtz equation for the smooth model in Fig.~\ref{fig:smooth_model}.}\label{table:test_smooth}
\end{table*}

\subsection{Fault model}

In general, iterative methods are very sensitive to discontinuous media.
At high frequency, interaction with short-wavelength structures, such as discontinuities, increases the number of reflections.
Each additional reflection requires additional iterations to convergence, hindering the efficiency of iterative methods.

\begin{figure}[h]
    \begin{center}
        \includegraphics[trim = 20mm 50mm 30mm 60mm, clip, width=8.5cm]{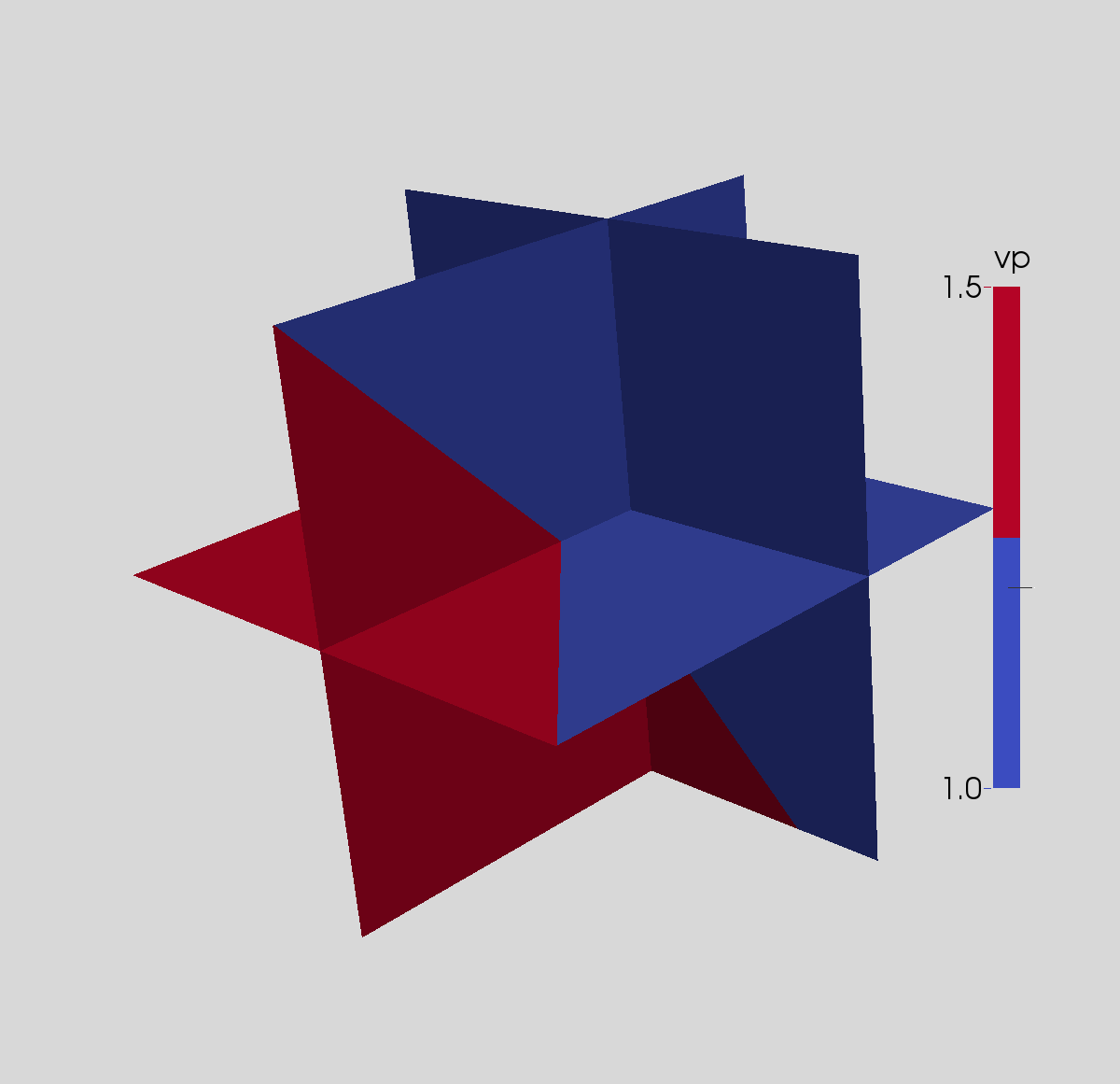}
    \end{center}
    \caption{Simple fault model.} \label{fig:fault_model}
\end{figure}

\begin{figure}[h]
    \begin{center}
        \includegraphics[trim = 20mm 50mm 30mm 60mm, clip, width=8.5cm]{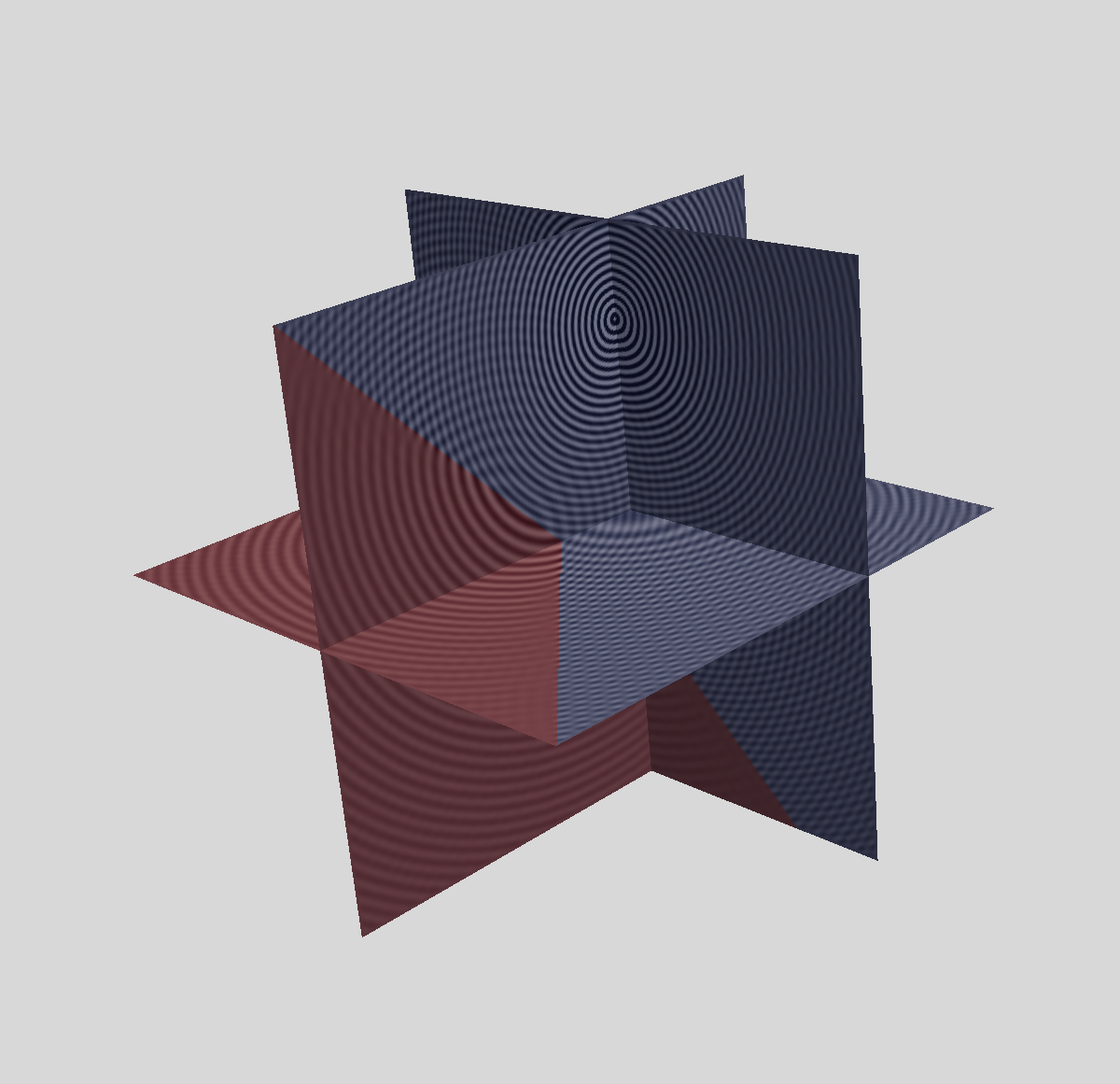}
    \end{center}
    \caption{Real part of the solution to the Helmholtz equation using the simple fault model.} \label{fig:solution_fault_model}
\end{figure}

Using the same configuration as for the homogeneous model, with the discontinuous velocity given in Fig.~\ref{fig:fault_model}, we demonstrate that the method of polarized traces deteriorates only marginally as a function of the frequency and number of subdomains.
A solution at 64 Hz is given in Fig.~\ref{fig:solution_fault_model} and the run-time scalability is again given in Figs.~\ref{fig:not_seam_1} and~\ref{fig:not_seam_R}.
As shown in Table \ref{table:test_fault}, we observe the same behavior as in the previous cases and that the strong reflection is handed efficiently by the transmission and polarizing conditions.

\begin{table*}
\begin{center}
\begin{tabular}{|r|cccccccc|}
\hline
N ~      & $50^3$ &  $100^3$ &  $100^3$ &  $200^3$ & $200^3$ &  $400^3$ & $400^3$  & $400^3$ \\
\hline
L  ~    & 5 & 10 & 10 & 20 & 20 & 40 & 40 & 40 \\
\hline
MPI Tasks   & 5 & 10 & 10 & 80 & 80 & 640 & 640 & 640\\
OpenMP Threads per Task     & 1 & 1 & 2 & 1 & 2 & 1 & 2 & 3\\
Total Cores & 5 & 10 & 20 & 80 & 160 & 640 & 1280 & 1920\\
Total Nodes & 1 & 1 & 2 & 5 & 10 & 80 & 80 & 128\\
\hline
\textbf{Single rhs} &  &  &  &  &  &  &  & \\
\hline
\# GMRES Iterations & 4 & 5 & 5 & 5 & 5 & 6 & 6 & 6\\
 Initialization [s]         & 0.4 & 1.1 & 1.0 & 7.3 & 4.7 & 20.4 & 20.3 & 21.0\\
 Factorization [s]      & 3.8 & 40.4 & 22.1 & 152.2 & 79.9 & 317.6 & 199.5 & 152.5\\
 Online [s]                 & 3.7 & 46.2 & 26.2 & 188.5 & 109.8 & 713.2 & 395.8 & 315.6\\
 Average GMRES [s]  & 0.8 & 8.1 & 4.6 & 33.0 & 19.2 & 106.2 & 58.7 & 46.5\\
 \hline
 \textbf{Pipelined rhs} &  &  &  &  &  &  &  & \\
\hline
 $R$ (number of rhs)    & 5 & 10 & 10 & 20 & 20 & 40 & 40 & 40\\
 Online [s]                     & 13.7 & 226.7 & 122.4 & 1222.7 & 647.1 & 4031.6 & 2710.6 & 1838.9\\
 Average GMRES [s]  & 2.9 & 40.1 & 21.6 & 216.5 & 114.7 & 605.0 & 409.9 & 276.3\\
 Online per rhs [s]                 & 2.7 & 22.7 & 12.2 & 61.1 & 32.4 & 100.8 & 67.7 & 46.0\\
 Average GMRES per rhs [s]  & 0.6 & 4.0 & 2.2 & 10.8 & 5.7 & 15.1 & 10.2 & 6.9 \\
 \hline
\end{tabular}
\end{center}
    \caption{Runtime (in seconds) for one and several right-hand sides for the solution of the Helmholtz equation for the fault model.}\label{table:test_fault}
\end{table*}

\subsection{SEAM model}
\begin{figure}[h]
    \begin{center}
        \includegraphics[trim = 20mm 90mm 30mm 80mm, clip, width=8.5cm]{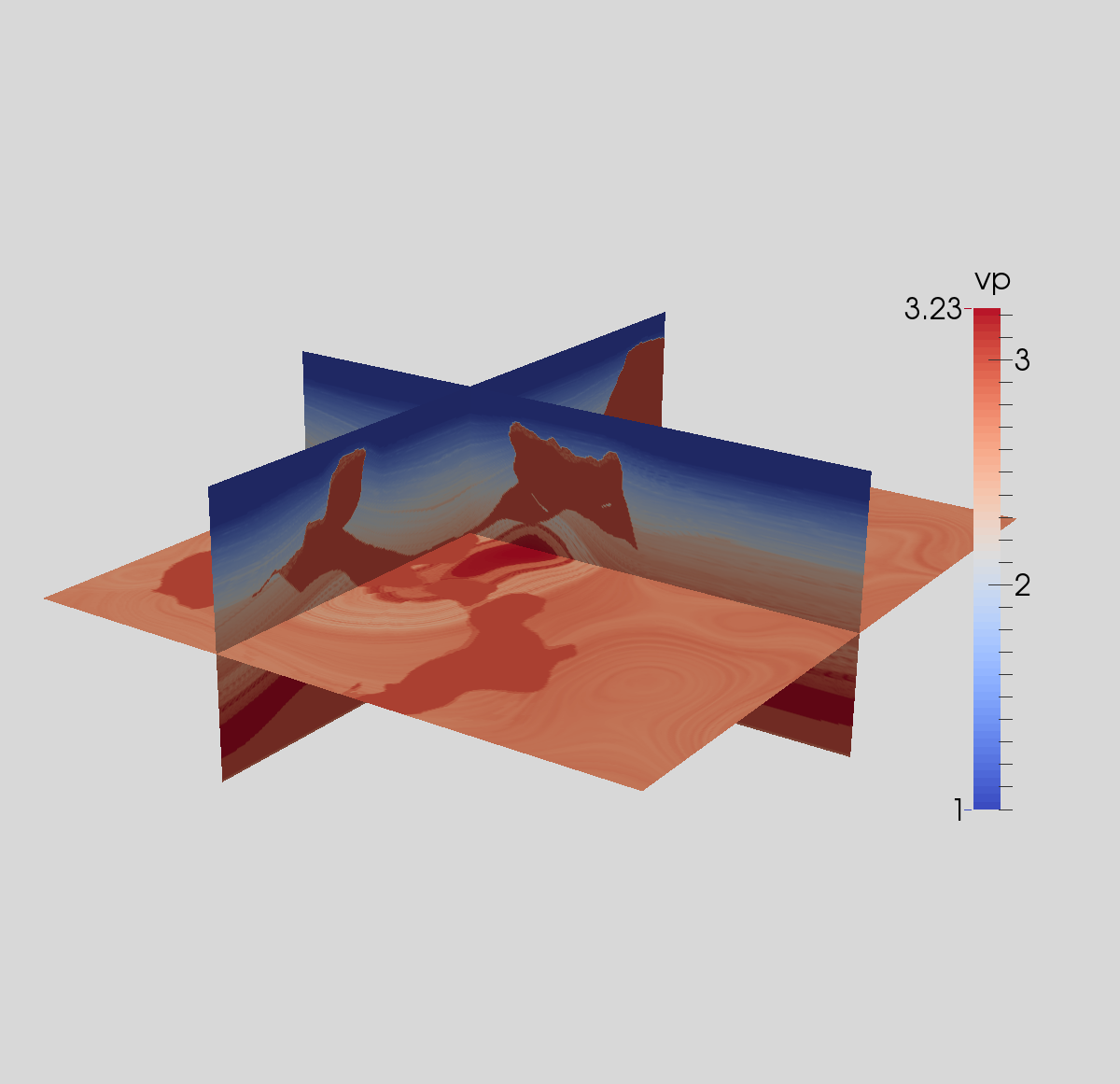}
    \end{center}
    \caption{Plot of the SEAM model.} \label{fig:SEAM_model}
\end{figure}

Beyond mere sensitivity to discontinuities of the medium, iterative solvers are highly sensitive to the roughness and heterogeneity of the velocity model, due to the great amount of interactions, reflections, and drastic changes of direction of waves due to high gradients in the wavespeed.
However, the performance method of polarized traces degrades only marginally for highly heterogeneous media, excepting resonant cavities.
We demonstrate this desirable performance on the SEAM Phase I velocity model (Fig.~\ref{fig:SEAM_model}).
In this experiment, we test 3 problem sizes, $N=0.65$M, $5.16$M, and $41.2$M degrees of freedom, which use $L=12, 24,$ and $48$ layers, respectively.
The remainder of the experimental setup is unchanged.

\begin{figure}[t]
    \begin{center}
        \includegraphics[clip,width=8.5cm]{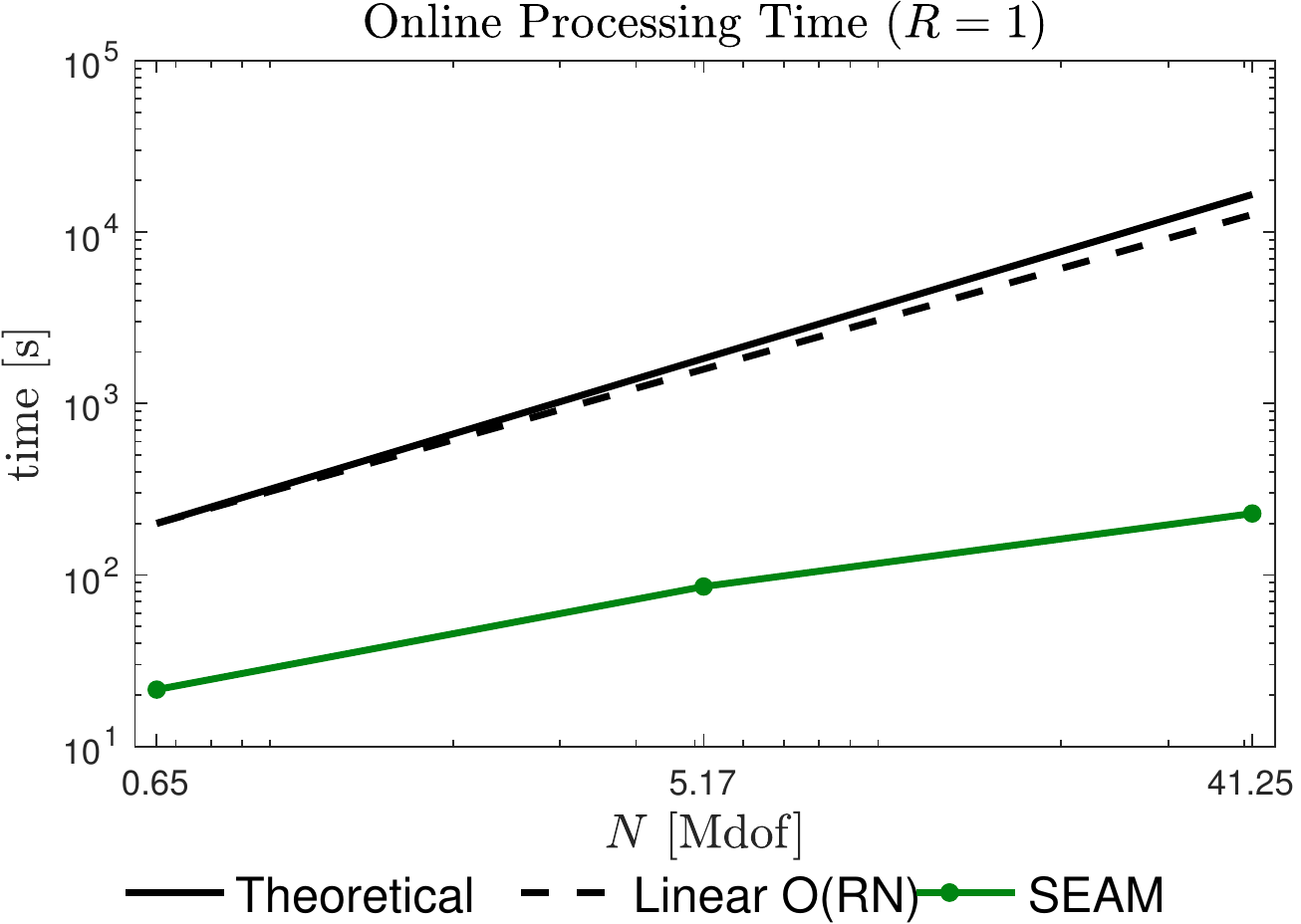}
    \end{center}
    \caption{Observed runtime as function of $N$ for the SEAM model for $R=1$ right-hand side.  For comparison, we show the theoretical scaling of the polarized traces algorithm (solid black), as well as a linear scaling (dashed black).} \label{fig:seam_1}
\end{figure}

\begin{figure}[t]
    \begin{center}
        \includegraphics[clip,width=8.5cm]{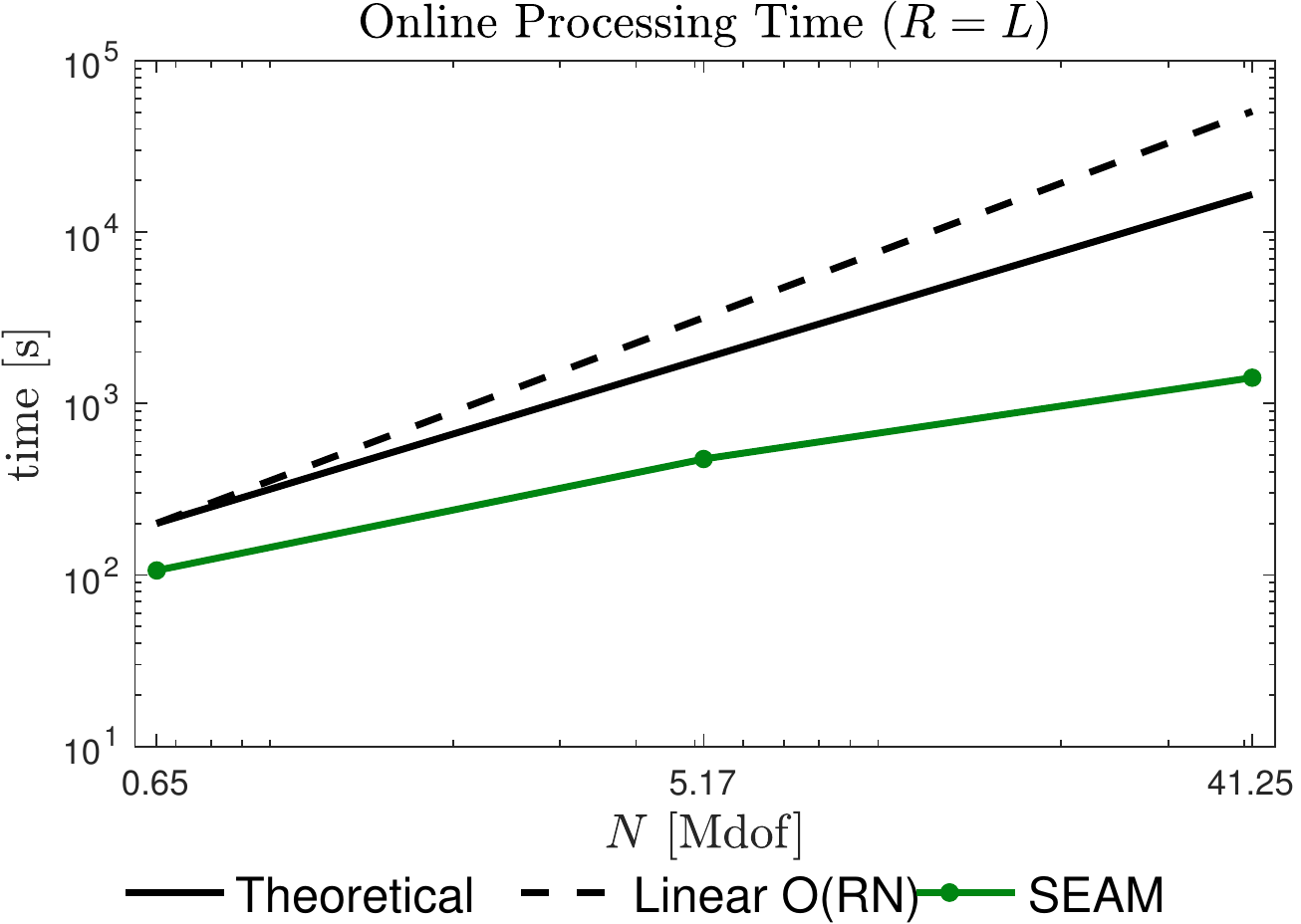}
    \end{center}
    \caption{Observed runtime as function of $N$ for the SEAM model for $R=L$ right-hand sides.  For comparison, we show the theoretical scaling of the polarized traces algorithm (solid black), as well as a linear scaling (dashed black).} \label{fig:seam_R}
\end{figure}

As seen in the data in Table \ref{table:test_SEAM} and plotted in Figs.~\ref{fig:seam_1} and~\ref{fig:seam_R}, the run-times are sub-linear with respect to the total number of unknowns.
Interestingly, we also observe a sub-linear run-time in the offline stages of the algorithm, which we attribute to the parallelism in the multi-frontal factorization.
This is expected, as the factorization is more computationally intensive than memory intensive.
In the more memory-intensive, and thus less parallel, solve phase, we still see an improvement over the theoretical curve, but the improvement is less pronounced.

\begin{figure}[t]
    \begin{center}
        \includegraphics[clip, width=8.5cm]{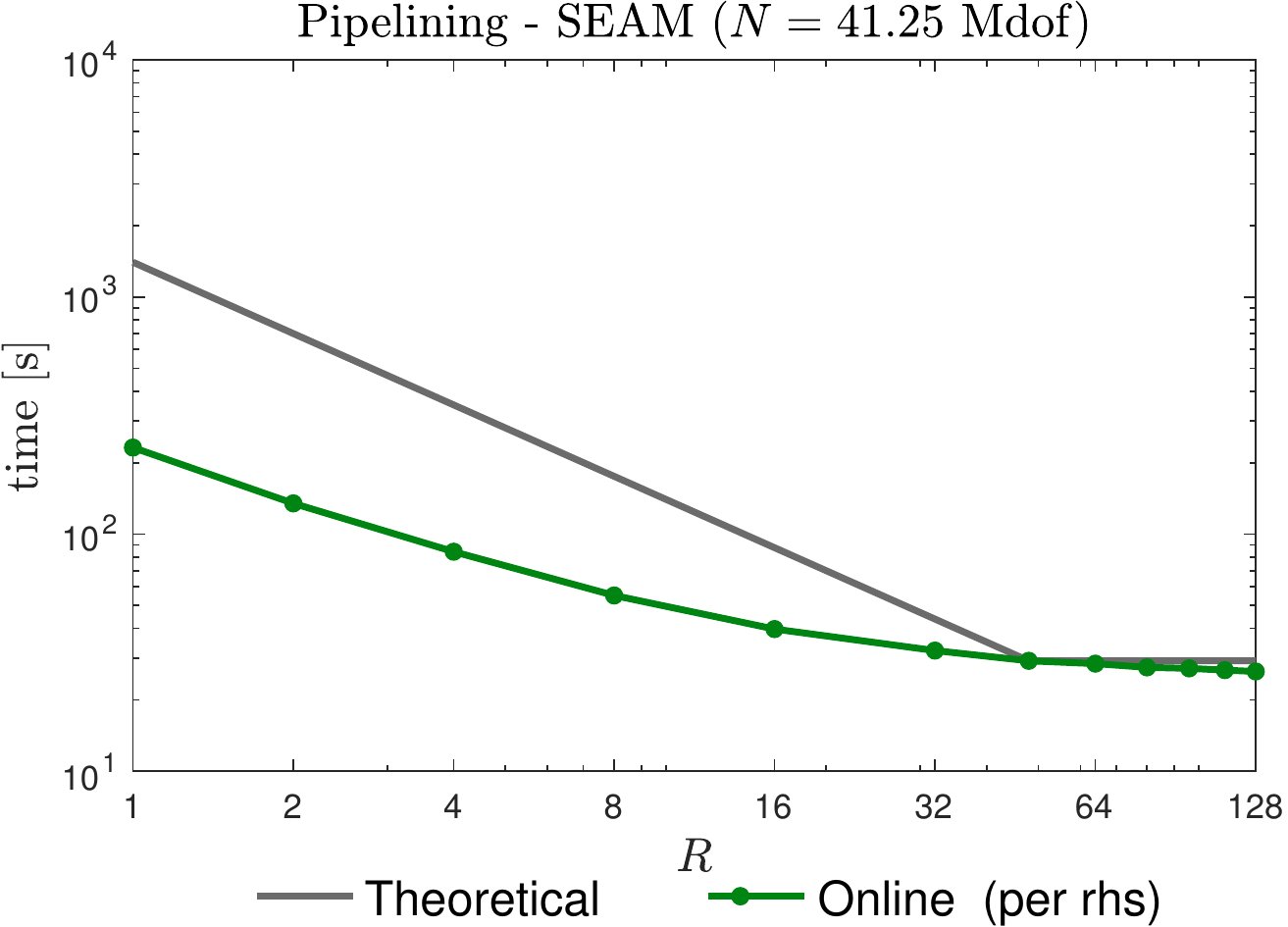}

    \end{center}
    \caption{Impact of scalability on pipelining for the SEAM model.  $N$ is held constant and $R$ is increased.  The observed runtime as function of $R$ is in green, and the dashed black line illustrates theoretical scalability.}
    \label{fig:pipeline_scaling}
\end{figure}

Finally, for the largest test case, we demonstrate the impact of pipelining by comparing the scalability of our method with the theoretical scalability, as a function of $R$.
As shown in Fig.~\ref{fig:pipeline_scaling}, experimental results indicate that we obtain the expected scalability.
Slight divergence from the theoretical curve is expected once the pipeline is fully saturated because the theoretical curve does not take into account the cost of filling and flushing the pipeline.

\begin{table*}[htb]
    \begin{center}
\begin{tabular}{|r|cccc|}
\hline
N ~     & $6.51\cdot 10^5 $  & $5.16\cdot 10^6$  & $4.12 \cdot 10^7 $ &  $4.12 \cdot 10^7 $  \\
\hline
L ~ & 12 & 24 & 48 & 48 \\
\hline
MPI Tasks       & 12 & 48 & 384 & 384\\
OpenMP Threads per Task     & 1 & 2 & 2 & 3\\
Total Cores     & 12 & 96 & 768 & 1152\\
Total Nodes     & 1 & 6 & 77 & 77\\
\hline
\textbf{Single rhs} &  &  &  & \\
\hline
 \# GMRES Iterations    & 4 & 5 & 6 & 6\\
 Initialization [s]     & 0.6 & 2.3 & 10.4 & 10.7\\
 Factorization [s]      & 15.2 & 46.5 & 111.4 & 97.9\\
 Online [s]             & 21.4 & 85.6 & 269.8 & 228.4\\
 Average GMRES [s]      & 4.6 & 14.9 & 40.0 & 33.7\\
 \hline
\textbf{Pipelined rhs} &  &  &  &  \\
\hline
 $R$ (number of rhs)            & 12 & 24 & 48 & 48\\
 Online [s]             & 106.3 & 474.8 & 1527.1 & 1415.4\\
 Average GMRES [s]  & 22.8 & 83.9 & 229.4 & 212.9\\
 Online per rhs [s]     & 8.8 & 19.8 & 31.8 & 29.5\\
 Average GMRES per rhs  [s] & 1.9 & 3.5 & 4.8 & 4.4\\
 \hline
\end{tabular}
\end{center}
    \caption{Runtime (in seconds) for one and several right-hand sides for the solution of the Helmholtz equation for the SEAM model.}\label{table:test_SEAM}
\end{table*}

\begin{figure}[h]
    \begin{center}
        \includegraphics[trim = 20mm 90mm 30mm 80mm, clip, width=8.5cm]{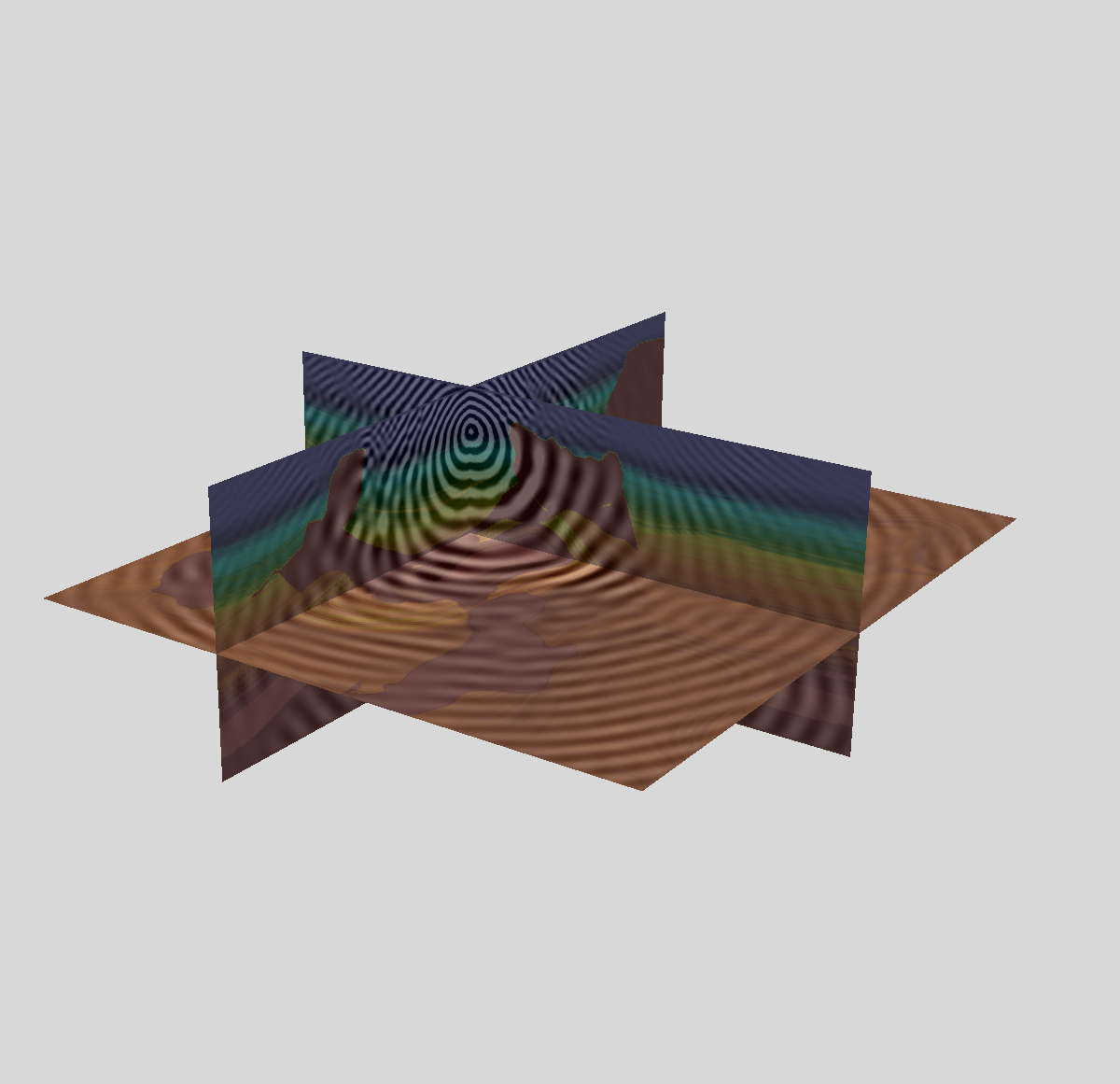}
    \end{center}
    \caption{Real part of the solution of the Helmholtz equation using the SEAM model.} \label{fig:solution_SEAM_model}
\end{figure}

Finally, Fig.~\ref{fig:pipeline_scaling} depics the behavior of the pipelining as we increase the number of right-hand sides. As expected, as we add more and more righ-hand sides to be solved simultaneously the runtime per right-hand side decreases, until the pipeline is full, when the average runtime remains almost constant.

\section{Conclusion}

We have presented a new and efficient solver for the $3$D high-frequency Helmholtz equation in heterogeneous media. The solver achieves a sub-linear runtimes by leveraging the solution of batches of right-hand sides properly pipelined and parallelism. The method presented in this paper broadens the applicability of parallel direct methods by embedding them in a domain decomposition method, whose rate of convergence is independent of the frequency.

Finally, the main limitations of the present method are large resonant cavities presenting high contrasts, for which the number of reflections can be large implying that the number of iterations for convergence can still be high.
\section{Acknowledgments}

This worked was sponsored by Total SA. LD is also sponsored by AFOSR grant FA9550-17-1-0316, ONR grant N00014-16-1-2122, and NSF grant DMS-1255203. The authors are grateful to Mike Fehler and Total SA for their permission to use the 3D SEAM model.

\bibliographystyle{seg}
\bibliography{GRF_integral_formulations.bib}

\appendix

\section{Appendix: Matrix-free algorithms} \label{appendix:matrix-free}

The application of the polarized system defined in Eq.~\ref{eq:polarized_SIE} is achieved by applying each block as shown in Alg.~\ref{alg:polarizedSystem}.

\begin{algorithm} Application of $ \mathbf{\underline{\underline{M}}} $ \label{alg:polarizedSystem}
    \begin{algorithmic}[1]
        \Function{ $\underline{\underline{\u}}$ = Application Polarized}{ $\underline{\underline{\v}}$ }
            \State $(\underline{\v}^{\downarrow},\underline{\v}^{\uparrow})  = \underline{\underline{\v}}$
            \State $\underline{\u}^{\downarrow} = \mathbf{\underline{D}}^{\downarrow} \v^{\downarrow}  + \mathbf{\underline{U}} \v^{\uparrow}$
            \State $\underline{\u}^{\uparrow} = \mathbf{\underline{D}}^{\uparrow} \v^{\uparrow}  + \mathbf{\underline{L}} \v^{\downarrow}$
            \State $\underline{\underline{\u}} = (\underline{\u}^{\downarrow},\underline{\u}^{\uparrow})$
        \EndFunction
    \end{algorithmic}
\end{algorithm}

To apply the blocks in a matrix-free fashion, we use Algs.~\ref{alg:downwards} and \ref{alg:upwardsReflections} in line 4 of Alg.~\ref{alg:polarizedSystem}, and  we use Algs.~\ref{alg:upwards} and \ref{alg:downwardsReflections} in  line 5 of Alg.~\ref{alg:polarizedSystem}. All algorithms in this section are embarrasingly parallel at the level of the layers as depicted in Fig.~\ref{fig:nodes_load}.

\begin{algorithm} Application of $ \mathbf{\underline{D}}^{\downarrow} $ \label{alg:downwards}
    \begin{algorithmic}[1]
        \Function{ $\underline{\u}^{\downarrow}$ = Downward Sweep}{ $\underline{\v}^{\downarrow}$ }
            \State $\u^{\downarrow,1}_{n^1}   = -\v^{\downarrow,1}_{n^1}$
            \State $\u^{\downarrow,1}_{n^1+1} = -\v^{\downarrow,1}_{n^1+1} $
            \For{  $\ell = 2: L-1$ }
                \State $ \tilde{\f}^{\ell} =  \delta(z_{1}-z)\v^{\downarrow,\ell-1}_{n^{\ell-1}} - \delta(z_{0}-z)\v^{\downarrow,\ell-1}_{n^{\ell-1}+1}  $
                \State $ \mathbf{w}^{\ell}  = (\mathbf{H}^{\ell})^{-1} \tilde{\f}^{\ell} $
                \State $\u^{\downarrow,\ell}_{n^{\ell}}   = \mathbf{w}_{n^{\ell}} - \v^{\downarrow,\ell}_{n^{\ell}}$
                \State $\u^{\downarrow,\ell}_{n^{\ell}+1} = \mathbf{w}_{n^{\ell}+1}-\v^{\downarrow,\ell}_{n^{\ell}+1} $
            \EndFor
            \State $\underline{\u}^{\downarrow} =  \left (\u^{\downarrow,1}_{n^1} , \u^{\downarrow,1}_{n^1+1}, \u^{\downarrow,2}_{n^2}, ..., \u^{\downarrow,L-1}_{n^{L-1}}, \u^{\downarrow,L-1}_{n^{L-1}+1}  \right)^{t} $
        \EndFunction
    \end{algorithmic}
\end{algorithm}

\begin{algorithm} Upward sweep, application of $ ( \mathbf{\underline{D}}^{\uparrow}  )^{-1}$  \label{alg:upwards}
    \begin{algorithmic}[1]
        \Function{ $\underline{\u}^{\uparrow}$ = Upward sweep}{ $\underline{\v}^{\uparrow}$ }
            \State $\u^{\uparrow,L}_{0}   = -\v^{\uparrow,L}_{0}$
            \State $\u^{\uparrow,L}_{1} = -\v^{\uparrow,L}_{1} $
            \For{  $\ell = L-1:2$ }
                \State $ \tilde{\f}^{\ell} = - \delta(z_{n^{\ell}+1}-z)\v^{\uparrow,\ell+1}_{0}  + \delta(z_{n^{\ell}}-z)\v^{\uparrow,\ell+1}_{1} $
                \State $ \mathbf{w}^{\ell}  = (\mathbf{H}^{\ell})^{-1} \tilde{\f}^{\ell}  $
                \State $\u^{\uparrow,\ell}_{1}   = \mathbf{w}^{\ell}_{1} - \v^{\uparrow,\ell}_{1}$
                \State $\u^{\uparrow,\ell}_{0} = \mathbf{w}^{\ell}_{0}-\v^{\uparrow,\ell}_{0} $
            \EndFor
            \State $\underline{\u}^{\uparrow} =  \left (\u^{\uparrow,2}_{0} , \u^{\uparrow,2}_{1}, \u^{\uparrow,3}_{0}, ..., \u^{\uparrow,L}_{0}, \u^{\uparrow,L}_{1}  \right)^{t} $
        \EndFunction
    \end{algorithmic}
\end{algorithm}

\begin{algorithm} Downwards reflections, application of $ \mathbf{\underline{U}}$ \label{alg:downwardsReflections}
    \begin{algorithmic}[1]
        \Function{ $\underline{\u}^{\uparrow}$ = Upward Reflections}{ $\underline{\v}^{\downarrow}$ }
            \For{  $\ell = 2:L-1$ }
                \State $ \begin{array}{ll} \f^{\ell} = & \delta(z_{1}-z)\v^{\downarrow,\ell}_{0} - \delta(z_{0}-z)\v^{\downarrow,\ell}_{1}  \\
                                                       & - \delta(z_{n^{\ell}+1}-z)\v^{\downarrow,\ell+1}_{1}  + \delta(z_{n^{\ell}}-z)\v^{\downarrow,\ell+1}_{0}
                                                    \end{array}$
                \State $ \mathbf{w}^{\ell}  = (\mathbf{H}^{\ell})^{-1} \f^{\ell} $
                \State $\u^{\uparrow,\ell}_{1} = \mathbf{w}^{\ell}_{1} - \v^{\downarrow,\ell}_{1}$
                \State $\u^{\uparrow,\ell}_{0} = \mathbf{w}^{\ell}_{0}  $
            \EndFor
            \State $ \f^{L} =  \delta(z_{1}-z)\v^{\uparrow,L}_{0} - \delta(z_{0}-z)\v^{\uparrow,L}_{1} $
                \State $ \mathbf{w}^{L}  = (\mathbf{H}^{L})^{-1} \f^{L} $
                \State $\u^{\uparrow,L}_{1} = \mathbf{w}^{L}_{1} - \v^{\downarrow,L}_{1}$
                \State $\u^{\uparrow,L}_{0} = \mathbf{w}^{L}_{0}  $
            \State $\underline{\u}^{\uparrow} =  \left (\u^{\uparrow,2}_{0} , \u^{\uparrow,2}_{1}, \u^{\uparrow,3}_{0}, ..., \u^{\uparrow,L}_{0}, \u^{\uparrow,L}_{1}  \right)^{t} $
        \EndFunction
    \end{algorithmic}
\end{algorithm}

\end{document}